\documentclass [hidelinks,12pt]{amsart}

\usepackage {amsmath,amssymb}
\usepackage {verbatim}
\usepackage {enumerate}
\usepackage {xcolor}
\usepackage {xspace}
\usepackage {todonotes}
\usepackage {mathrsfs}
\usepackage [all]{xy}
\usepackage {hyperref}
\hypersetup {
    colorlinks,
    linkcolor={blue!70!black},
    citecolor={blue!70!black},
    urlcolor={teal!50!black}
}

\usepackage [normalem]{ulem} % for the command \sout {}, which strikes out the text inside " {  }"

\newcommand {\mult }{\mathfrak {M}}
\newcommand {\munit }{{e}}
\newcommand {\prp }{partial representation property\xspace }
\newcommand {\m }{{}^{-1}}
\newcommand {\ps }[1]{I^{#1}}

\newcommand {\fac }[1]{\textcolor {blue!75!black}{#1}}%fernando
\theoremstyle {plain}
\newtheorem {thm}{Theorem}[section]
\newtheorem {prop}[thm]{Proposition}

\newtheorem {cor}[thm]{Corollary}

\newtheorem {thm*}{Theorem}
\newtheorem {prop*}{Proposition}
\newtheorem {lem*}{Lemma}
\newtheorem {cor*}{Corollary}

\theoremstyle {definition}
\newtheorem {defn}[thm]{Definition}

\theoremstyle {remark}
\newtheorem {rem}[thm]{Remark}

\title [Strong equivalence of graded algebras]
{Strong equivalence of graded algebras}
\author {F.~Abadie}
\address{Centro de Matem{\'a}tica, Facultad de Ciencias,
 Universidad de la Rep{\'u}blica, Igu\'a 4225, CP 11400,
 Montevideo, Uruguay}
\email{fabadie@cmat.edu.uy}
\author {M.~Dokuchaev}
\address{Departamento de Matem\'atica, Universidade de S\~ao Paulo,
Rua do Mat\~ao, 1010, 05508-090 S\~ao Paulo, Brazil}
\email{dokucha@ime.usp.br}
\author {R.~Exel}
\address{Departamento de Matem\'atica, Universidade Federal de Santa Catarina, 88040-900 Florian{\'o}polis, Brazil}
\email{ruyexel@gmail.com}

\begin {document}
%\begin {comment}

%%%%%%%%%%%%%%%%%%%%%%%%%%%%%%%%%%%%%%%%%%%%%%%%%%%%%%%%%%%%%%%%%%%%%%%%%%%%%%%
%% Comentarios de Ruy - Dezembro de 2021

% \vspace *{-3cm}
% \begin {flushright}
% \footnotesize {\textsf {Preliminary version, \today }}
% \end {flushright}
% \vspace *{3cm}

\begin {abstract}
We introduce the notion of a strong equivalence  between
graded algebras and prove that any partially-strongly-graded algebra
by a group $G$
is strongly-graded-equi\-va\-lent
to the skew group algebra by a product partial action of $G.$
 %and graded equivalent to a global skew group algebra.
%As a consequence, any strongly graded algebra is strongly-graded-equivalent to a skew group algebra.
  As to a more  ge\-ne\-ral idempotent graded algebra $B,$ we point out that
the Cohen-Montgomery duality holds for $B$, and $B$ is
 graded-equi\-va\-lent to a global skew group algebra.
%(the concept of the graded equi\-va\-len\-ce is adapted  for idempotent graded algebras in terms of  a Morita context).
We show that  strongly-graded-equi\-va\-len\-ce preserves strong
  gradings and is nicely related to Morita equivalence of pro\-duct partial actions.
  Furthermore,
we prove that any product partial group action  $\alpha $ is globalizable up to Morita equivalence; if such a
glo\-ba\-li\-za\-tion $\beta $ is minimal, then the skew group algebras by $\alpha $ and $\beta $ are
graded-equivalent;
moreover, $\beta $ is unique up to Morita equivalence. Finally, we show that strongly-graded-equivalent
 partially-strongly-graded algebras {\color{black} with ortho\-go\-nal local units} are stably isomorphic as graded algebras.
\end {abstract}

%\end {comment}
\maketitle
{\color {blue!70!black}
\tableofcontents }
\section {Introduction}

Two kinds of graded-equivalences are usually  considered when dealing with graded rings: one, the graded-equivalence,
as defined  in \cite {GordonGreen} and \cite {MeniniNast}, and a stronger one, the graded Morita equivalence, given in
\cite {Boisen}.  A general description of equivalences of categories of graded modules over unital rings graded by
different groups was given in \cite {ADelRio91}. Morita theory for unital rings was extended to rings with local units in \cite {Abrams} and \cite {AnhMarki}, with  graded versions worked out in  \cite {Haefner94} (graded-equivalence) and
\cite {Haefner95} (graded Morita equivalence). To a great extent the graded theory is stimulated by the
Cohen-Montgomery duality, saying that if $B$ is a unital ring graded by a finite group $G$ of order $n,$ then
the skew group ring     $(B \#G)\rtimes _\beta G$ and the matrix ring $M_n(B)$ are isomorphic,
where $\beta $ is a certain canonical global action of $G$ on
the smash product  $B \#G .$ Extensions of the above duality theorem were obtained, in particular,  in \cite {quinn}, \cite {beattiedual}, \cite {AlbuNast},  \cite {ADelRio92} and  \cite {Haefner94}.

 %In \cite [Corollary 2.10]{Haefner} J. Haefner proved that a graded ring with graded local units is graded-equivalent to a skew %group ring. We give the following.

Cohen-Montgomery duality was inspired by the use of duality  in studying von Neumann algebras and C*-algebras.
 In C*-theory a duality development based on partial group actions was done in \cite {abenv}, introducing,
in particular,
the notion of the Morita equivalence of partial actions and stimulating algebraic analogues in \cite{ades}.
  Partial actions on algebras
 are closely related with graded algebras: on the one hand any partial group action on an algebra gives rise to the corresponding crossed product, and on the other by \cite{des} an algebra partially strongly graded  by an arbitrary group $G$ (see Definition~\ref {defn:prp} below) with enough local units is stably isomorphic to a crossed product by a twisted partial action of $G.$ The latter fact is a purely  algebraic analogue of a C*-algebraic result from \cite{extwist} stated in terms of an important concept of a Fell bundle, which is roughly the essential data remaining after disassembling a graded C*-algebra into its homogeneous pieces.
Moreover, any Fell bundle carries an associated
   partial action of the underlying group on the spectrum of the unit
    fiber, as shown in \cite {aa}.
 % the so
 %   called Fell bundles over groups, which
    %in turn can be thougt of
    %as an abstraction of gradings of C*-algebras. On the one hand, any
    %partial action induces a Fell bundle, known as the semidirect
    %product bundle. On the other hand, it was shown in \cite {extwist}
    %that many Fell bundles can be described as semidirect product
    %bundles of twisted partial actions, while the twist was later shown to be %superfluous in \cite {sehnem}
    %(see also \cite [Chapter 27]{pdsfba}).
    %The strong
    %relationship between Fell bundles and partial actions can be
    %illustrated by the fact that any Fell bundle carries an associated
 %   partial action of the underlying group on the spectrum of the unit
 %   fiber, as shown in \cite {aa}.

 More recently, notions of \textit {weak}
    equivalence and \textit {strong} equivalence between Fell
    bundles were introduced in \cite {af} and \cite {abf},
    %where it was
    %shown that any Fell bundle is strongly-equivalent to the
   % semidirect bundle of a partial action, and weakly-equivalent to
   % the semidirect Fell bundle of a global action. This
    %latter fact should be interpreted as a globalization result
   % because, at least heuristically, a Fell bundle can be seen as a
  %  kind of partial action by Hilbert bimodules \cite {BussMeyerZhu} (from this point
   % of view the notion of strong equivalence corresponds to an
   % equivariant morphism).
and in the present work we look at these notions from a purely algebraic point of view.
%versions of  the above results.
In the algebraic case, a Fell bundle
    is essentially a graded algebra, more precisely a
    \textit {partially-strongly-graded algebra}, as
    introduced in  Definition~\ref {defn:prp} below, and weak
    equivalence  corresponds to graded-equivalence. However, strong
    equivalence gives rise to a new kind of relation between partially-strongly-graded algebras, which we call here strong-graded-equivalence.

    Crossed products by partial actions
   on algebras came into algebra from the theory of C*-algebras. Since a C*-algebra enjoys the nice property that the
   intersection of two closed ideals is equal to their product,
   the notion of a partial action in this category can be
   introduced equivalently in terms of intersection or product of
   ideals. However the situation is very different in the purely
   algebraic framework, and one has to make a choice.
   This fact was not explicitly observed in the first
   algebraic works on partial actions, but it emerged clearly when
   twisted partial actions were considered in \cite {des}. It appeared again when the question of
   glo\-ba\-li\-za\-tion up to Morita equivalence was considered in
   \cite {ades}, giving rise to the notion of a \textit {regular partial
     action}. The work done in the present article suggests that
   perhaps the notion of a product partial action is the most convenient
   choice to translate the concept of a partial action from the
   C*-algebraic world to the purely algebraic one.

 We begin the paper recalling some notions on graded algebras in Section~\ref {Sec:GrGradedAlgs}
and pointing out in  Theorem~\ref {thm:duality} that
the above mentioned duality   holds for any idempotent  algebra $B$ graded by an arbitrary  group $G,$ with $M_n(B)$
 replaced
by the algebra ${\rm FMat}_{G} (B)$ of $G\times G$-matrices over $B$ with only a finite number of non-zero entries.
For graded algebras with $1$ this is known by  \cite [Theorem~2.2]{beattiedual}. In Section~\ref {sec:GrEquiv} we define graded-equivalence of
idempotent graded algebras using graded Morita contexts and prove in Theorem~\ref {thm:geq} that if $B$ is an
idempotent algebra  graded by a group $G$, then $B$ and the skew group algebra  $(B\# G)\rtimes _{{\beta }^B}G$
are graded-equivalent,
where $B\# G$ is the  Beattie's version of the smash product \cite {beattiedual} and $ {\beta }^B$ is the usual global
action of $G$ on $B\# G$, sometimes referred to as the {\it dual} action.

 Strong-graded-equivalence is introduced in
 Section~\ref {sec:StrGrEquiv}, and we deal with it concentrating  on
 {\it partially-strongly-graded algebras} $A=\oplus _{t\in G} A_t,$ i.e. we assume that   the equality $A_t=A_tA_{t^{-1}}A_t$ is satisfied  for all
$t\in G.$ The latter condition naturally appeared in a characterization of graded algebras as crossed products  by twisted partial group actions in \cite {des}. Section~\ref {sec:ppa} is dedicated to product partial actions, the main result being
Theorem~\ref {thm:crossed product ppa and env}, which says that the skew group algebra $A\rtimes _\alpha G$
by a product partial action $\alpha $ is graded-equivalent to the skew group algebra
$B\rtimes _\beta G,$ where $\beta $  is a minimal globalization  of $\alpha .$
Given a  partially-strongly-graded algebra  $B=\oplus _{t\in G}B_t,$ the dual action $\beta ^B$ of $G$ on $B\#G$ can
be restricted to the ideal  $\ps {B}$ of  $B\#G$ defined in
Section~\ref {Sec:GrGradedAlgs} (which is called \textit {partial smash
  product} in this work), resulting in a  pro\-duct partial action
$\gamma ^B,$ called the {\it canonical partial action} associated to
$B$, such that $B$ and the skew group algebra $I^B\rtimes _{\gamma ^B}G$ are
strongly-graded-equivalent.
Theorem~\ref {thm:partialrep} states that  the (global) skew group algebra
$(B\#G) \rtimes _{\beta ^B} G$ and the partial one $\ps {B}\rtimes _{{\gamma }^B} G$ are graded-equivalent. We begin Section~\ref {MoritaEqParAc}   by adapting to product partial actions the Morita equivalence of  partial actions
considered earlier in \cite {ades}. One of the main results of the section is Theorem~\ref {thm:sg}, which  says that
if $B=\oplus _{t\in G}B_t$ is a partially-strongly-graded algebra, then $B$
is strongly-graded-equivalent to $\ps {B}\rtimes _{{\gamma }^B} G$.  As a consequence we obtain in
Corollary~\ref {cor:sg} that
if $B$ is a strongly-graded algebra, then $B$ is
strongly-graded-equivalent to the (global) skew group algebra  $(B\# G) \rtimes _{\beta ^B}G.$
 Another consequence states that the crossed product by any twisted partial group action is strongly-graded-equivalent to the skew group algebra of a product partial
 action (see Corollary~\ref {cor:tcp}).  We also prove in   Theorem~\ref {thm:strongmoritaE} that
the canonical partial actions $\gamma ^A$ and $\gamma ^B$ associated to strongly-graded-equivalent  partially-strongly-$G$-graded
algebras $A$ and
$B$, are  Morita equivalent. We also consider, at the end of
Section~\ref {MoritaEqParAc},
  another notion of equivalence of product partial actions,
  weaker than Morita equivalence.

 The question of whether a partial action is the restriction of
   a global action, that is, the question of \textit {globalization},
   is one of the most important topics in the theory of partial
   actions. It was initially considered in
   \cite {abenv} in the context of C*-algebras, and  afterwards   in a
   series of algebraic  papers.
In  Section~\ref {sec:glob} we deal with globalization up to Morita equivalence. More specifically, we show
in Theorem~\ref {thm:main2} that
 any pro\-duct partial action $\alpha $ of a group $G$ has a so-called Morita enveloping action
$\beta ,$  i.e. $\beta $ is a minimal globalization of  a product partial action which  is Morita equivalent to  $\alpha .$
Furthermore, $\beta $
  is unique up to Morita equivalence, and the skew group algebras $A \rtimes _\alpha
G$ and $B \rtimes _\beta G$ are graded-equivalent. In order to  prove Theorem~\ref {thm:main2} we
establish several  facts, one of them being Theorem~\ref {thm:moritaskew}, saying  that  skew group algebras by
pro\-duct partial actions $\alpha $ and $\alpha '$ are strongly-equivalent if and only if  $\alpha $ and $\alpha '$ are
 Morita equivalent. We then use our results on globalization to give
 different characterizations of Morita equivalence and weak
 equivalence of product partial actions. Finally, in Section~\ref {sec:stab} we  employ the technique developed in
\cite {des} to prove in
Theorem~\ref {thm:gradedstableiso} that given strongly-graded-equivalent partially-strongly-$G$-graded algebras
$A$ and $B$
with orthogonal local units,    there exists a graded isomorphism of algebras
${\rm FMat}_{\mathcal X } (A) \cong {\rm FMat}_{\mathcal X} (B),$ where ${\mathcal X}$ is a sufficiently large cardinal.

 In what follows $G$ will stand for an arbitrary group and  $k$ for  an arbitrary commutative
associative
unital ring, which will be the base ring for our algebras. The latter
will be assumed to be associative  and non-necessarily unital.
Let $A$ and $B$ be algebras. A left module $_A M$ over  $A$ is said to
be {\it unital} %\todo {the   term ``firm module'' is excluded as being
                %a stronger condition than that of a unital module}
if $AM =M.$ We shall say that an $(A,B)$-bimodule $_A M _B$ is {\it unital} if $AM =M = M B.$
Given  a right $A$-module $M_A,$  a left $A$-module $_A N$ and subsets  $M'\subseteq M, N' \subseteq N,$
we denote by
$M' \otimes _A N'$ the $k$-submodule of the  tensor product
$M\otimes _A N$  generated by all elements of the form $x \otimes y $ with $x \in M', y \in N'.$

\section {Group graded algebras}\label {Sec:GrGradedAlgs}
Let $G$ be a group, and let $B=\oplus _{t\in G}B_t$ be a $G$-graded
algebra, possibly non-unital. We denote by $1$ the
unit element of $G.$ If  $b\in B,$ we write $b_t$ for the homogeneous component of $b$ in $B_t,$
so that $b = \sum _{t\in G} b_t .$ Note also that $b_t$ may also stand for an element in $B_t$ not necessarily
related to some $b\in B.$ This will be clear from the context and no confusion should arise. Consider the algebra
${\rm RCFMat}_G(B)%=\{(b_{r,s})_{r,s\in G}: b_{r,s}\in B\}
$ of all row and column finite  $G\times G$-matrices with coefficients in $B$, with the usual operations of addition and
multiplication of matrices.
%\sout{This is a
%  $G$-graded algebra: ${\rm RCFMat}_G(B)=\oplus _{t\in G}R_t$, where
%  $R_t:=\{M\in {\rm RCFMat}_G(B): M(r,s)=0 \textrm { if }rs^{-1}\neq
%  t\}$.}
An $(r, s)$-position of a matrix $d\in \rm{RCFMat}_G(B)$ will be denoted by $d(r, s)$. Observe that the algebra ${\rm RCFMat}_G(B)$  is
unital if so is $B$.  Let ${\rm FMat}_G (B)$ be the two-sided
ideal of ${\rm RCFMat}_G(B)$  whose elements are the matrices with finitely many non-zero entries.   If $B$ is unital then $\munit _{r,s} \in {\rm RCFMat}_G(B)$ will denote the matrix unit and   $b \munit _{r,s},$ with $b\in B,$ will stand for the product of the scalar matrix $bI$ with  $\munit _{r,s},$  where $I\in {\rm RCFMat}_G(B)$ is the identity matrix. Then clearly      $b \munit _{r,s}$ is the matrix having $b$ at $(r,s)$-position and $0$ at all other positions. By convention we shall denote this matrix by $b \munit _{r,s}$ even if $B$ is non-unital.
% We denote by $b_{r,s} \munit _{r,s}$ the matrix whose only
%possible non-zero entry is precisely $b_{r,s}\in B$, at the
%$(r,s)$-position.
% Let $C$ be any unital algebra that contains $B$ as a two-sided ideal. In this case ${\rm RCFMat}_G(C)$ is a
%$B$-bimodule in an obvious way, and $b_{r,s}\munit _{r,s}$ can be seen as the element $\munit _{r,s}\in
%$={\rm RCFMat}_G(C)$ acted on the left by $b_{r,s}$.
 Then ${\rm FMat}_G(B)=
 \textrm {span}\{b\munit _{r,s}: r,s\in G, b\in
B\}$. We will be interested in the following subalgebra $B\# G$ of
${\rm FMat}_G(B)$:
\[ B\# G:=\textrm {span}\{b \munit _{r,s}\in
{\rm FMat}_G(B): r,s\in G, b \in B_{r^{-1}s}\}.\]
Thus $B\# G=\oplus _{r,s\in G}B_{r^{-1}s}\munit _{r,s}$. Note that if $B$ is unital, then
this algebra is nothing but the smash product  in the sense of
Beattie~\cite {beattiedual}  (see also
\cite {quinn}),
  which in the
case of a finite  $G$  agrees with the smash product  in
\cite {cm}. If $A=\oplus _{t\in G}A_t$ is
another $G$-graded algebra, and $\phi :A\to B$ is a (graded) homomorphism of
$G$-graded algebras, then $\phi $ induces a homomorphism $\phi ^\#: A\#
G\to B\# G$ such that
$\phi ^\#(a_{r^{-1}s}\munit _{r,s})=\phi (a_{r^{-1}s})\munit _{r,s}$. It
is easily shown that $B\mapsto B\# G$, $\phi \mapsto \phi ^\#$ is a
functor from the category of $G$-graded algebras into the category of
algebras. Note that $\phi ^\#$ is injective if and only if so is
$\phi $. Moreover, suppose that $B=\oplus _{t\in G}B_t$ is a $G$-graded
algebra, and that $A=\oplus _{t\in G}A_t$ is a $G$-graded subalgebra of
$B$. Then, if $\iota :A\to B$ is the natural inclusion, we have that
$\iota ^\#:A\# G\to B\# G$ is also the natural inclusion.
\par
%{\color{blue} \sout{Similarly, if $N=\oplus _{t\in G}N_t$ is a $G$-graded module
%over  $B=\oplus _{t\in G}B_t$, we define ${\rm RCFMat}_G(N)$,
%${\rm FMat}_G(N)$ and $N\# G$ in the same way as done for $G$-graded
%algebras.   For instance,
%$N\# G:=\textrm {span}\{n_{r^{-1}s}\munit _{r,s}\in
%{\rm FMat}_G(N):\, n_{r^{-1}s}\in N_{r^{-1}s}\}$.
%Note that $N\#
%G$ is naturally a $B\# G$-module. In particular, in case $N$ is a
%$G$-graded left (or right) ideal in  $B$, then $N\# G$ %is also a left
%(respectively: right) ideal in $B\# G$.}}
%\par
We may think of
$B$ as a subalgebra of ${\rm RCFMat}_G(B)$ via the map $\eta :b\mapsto \eta (b)$
such that $\eta (b)(r,s)=b_{r^{-1}s}$, $\forall r, s\in G$.
If $G$ is
finite, then the map $\eta $ has its range contained in $B\# G$.
\par
There is a natural action  %\todo {The symbol $\beta $ was overused in
%the preprint, so in several places I replaced $\beta $ by $\beta ^B$,
%which is the restriction of this $\beta $ to the smash
%product. \textcolor {blue}{OK!} }
$\beta $ of $G$ on ${\rm RCFMat}_G(B)$, such that
$(t\cdot d)(r,s)=d(t^{-1}r,t^{-1}s)$, $\forall r,s,t\in G$, $d\in
{\rm RCFMat}_G(B)$. Thus  $t\cdot
(b \munit _{r,s})=b \munit _{tr,ts}$, $\forall
t\in G$, $b \munit _{r,s}\in {\rm RCFMat}_G(B)$. Clearly, the subalgebras
${\rm FMat}_G(B)$, $B\# G$, and $\eta (B)$ are
invariant under~$\beta $.
%In fact we have ${\rm RCFMat}_G(B)^{\beta }=\eta (B)$, that
%is, $\eta (B)$ is precisely the subalgebra of fixed points of $\beta $. '
We
denote by $\beta ^B$ the {\it dual action},  namely restriction of $\beta$  to
$B\# G$.  This action is natural with respect to the smash
product functor:
$\phi ^\#\beta _t^A=\beta _t^B\phi ^\#$, $\forall t\in G$.  Note that each
element of  $\eta (B)$ is fixed by $\beta $ and, moreover, if $G$ is finite, then
$\eta (B) \subseteq B\# G$ is precisely the subalgebra of fixed points of $\beta ^B$.

\par
We concentrate now our attention on the smash product
  $B\# G$ of
the $G$-graded algebra $B$.   The skew group algebra
$$(B\#G)\rtimes _{\color {black} {\beta }^B} G
= \oplus _{t \in G} (B\#G) \delta _t$$  possesses  the  usual  $G$-grading, defined by declaring  $(B\#G) \delta _t$ to be the $t$-homogeneous component of  $(B\#G)\rtimes _{\color {black} {\beta }^B} G,$ where $t\in G.$ In Theorem~\ref{thm:duality} below  we shall consider a $G$-grading on ${\rm FMat}_G(B),$ defined  by setting $R_u = {\rm span} \{ B_{r} \munit _{s,t} \,: \,  s r t\m = u \} $ to be the $u $-homogeneous component of ${\rm FMat}_G(B),$ where $u, r,s, t\in G.$
 It is readily verified that this indeed defines a $G$-grading.

Our first remark about $B\# G$ is
the following duality theorem, that generalizes \cite [Theorem~3.5]{cm}
and \cite [Theorem~2.2]{beattiedual} (see also
\cite [Theorem~1.3]{quinn}):
\begin {thm}\label {thm:duality}
  Let $B$ be any $G$-graded algebra. Then the skew group algebra
$(B\#G)\rtimes _{\color {black} {\beta }^B} G$  is
 naturally isomorphic, as a graded algebra,  to ${\rm FMat}_G(B).$
\end {thm}
\begin {proof}
Let $C$ be a unital algebra that contains $B$ as a two-sided
ideal. Thus ${\rm RCFMat}_G(B)$ and
${\rm FMat}_G(B)$ are two-sided ideals of ${\rm RCFMat}_G(C)$. For $t\in
G$, let $\Delta _t\in {\rm RCFMat}_G(C)$ be such that
$\Delta _t(r,s)=[r=ts]\in C,$ where the (square) brackets stand for the boolean value.
Let $\psi _B:(B\# G)\rtimes _{\color {black} {\beta }^B}
G\to {\rm FMat}_G(B)$ be the $(B\# G)$-module map given by
$\psi _B(c\delta _t)=c\Delta _t\in (B\# G){\rm RCFMat}_G(C)\subseteq
{\rm FMat}_G(B)$. We have  \[  (\munit _{r,s}\Delta _t )  (u,v)=\sum _{w\in
  G}\munit _{r,s}(u,w)\Delta _t(w,v)
=[r=u][s=tv]=\munit _{r,t^{-1}s}(u,v).\] Then
$\munit _{r,s}\Delta _t=\munit _{r,t^{-1}s}$ (similarly one can show that $\Delta _t\munit _{r,s}=\munit _{tr,s}$, so
$\beta _t( \munit _{r,s} ) = \Delta _t \munit _{r,s} \Delta _{t^{-1}} $).
 Now if $c_1=b_1 \munit _{r_1,s_1}$,
$c_2=b_2 \munit _{r_2,s_2}$,
$b_i \in B_{r_i^{-1}s_i} ,$ $ r_i, s_i , r,s\in G, i=1,2,$ then
\begin {align*}
\psi _B(c_1\delta _r)\psi _B(c_2\delta _s)%=(c_1\Delta _r)(c_2\Delta _s)
&=b_1 \munit _{r_1,s_1}  \Delta _r \,
               b_2 \munit _{r_2,s_2} \Delta _s\\
&=b_1 b_2 \munit _{r_1,r^{-1}s_1}\munit _{r_2,s^{-1}s_2}\\
&=[s_1=rr_2]b_1 b_2 \munit _{r_1,s^{-1}s_2}.
\end {align*}
On the other hand:
\begin {align*}
\psi _B((c_1\delta _r)(c_2\delta _s))
&=\psi _B((b_1 \munit _{r_1,s_1}\delta _r)(
               b_2 \munit _{r_2,s_2}\delta _s))\\
&=\psi _B(b_1 b_2 \munit _{r_1,s_1}\munit _{rr_2,rs_2}\delta _{rs})\\
&b_1 b_2 \munit _{r_1,s_1}\munit _{rr_2,rs_2}\Delta _{rs}\\
&b_1 b_2 \munit _{r_1,s_1}\munit _{rr_2,s^{-1}r^{-1}rs_2}\\
&[s_1=rr_2]b_1 b_2 \munit _{r_1,s^{-1}s_2}.
\end {align*}
 Hence $\psi _B$ is a homomorphism of algebras. Moreover, it
  respects the gradings, since
  $\psi _B (B_{r\m s} \munit _{r,s} \delta _t) =
  B_{r\m s} \munit _{r,s} \Delta _t = B_{r\m s} \munit _{r,t\m s},$ whose degree is $r (r\m  s) (t\m s)\m = t.$ Furthermore, $\psi _B$  is injective, because
$\{\Delta _t:t\in G\}$ is $(B\# G)$-linearly independent: if
$\sum _{t\in G}c_t\Delta _t$ is a finite sum which is equal to zero,
then $\sum _{t\in G}c_t\Delta _t(u,v)=0$, $\forall u,v\in G$, that is,
$\sum _{t\in G}c_t[u=tv]=0$, $\forall u,v\in G$; fixing $t$ and choosing
$u=t$, $v=e$, we conclude $c_t=0$. Let us
compute now the range of $\psi _B$:
\begin {gather*}
\psi _B((B\# G)\rtimes _{\color {black} {\beta }^B} G)
=\sum _{t\in G}(B\# G)\Delta _t
=\sum _{t\in G}\oplus _{r,s\in G}B_{r^{-1}s}\munit _{r,s}\Delta _t\\
=\sum _{t\in G}\sum _{r,s\in G}B_{r^{-1}s}\munit _{r,t^{-1}s}
=\oplus _{r,u\in G}\sum _{t\in G}B_{r^{-1}tu}\munit _{r,u}\\
=\oplus _{r,s\in G}\sum _{v\in G}B_v\munit _{r,s}
=\oplus _{r,s}B\munit _{r,s} ={\rm FMat}_G(B).
\end {gather*} This ends the proof
of $( B\# G ) \rtimes _{\color {black} {\beta }^B} G\cong {\rm FMat}_G(B)$.
\par
We next show the naturality of the isomorphism $\psi _B$.  The definition of
 $\psi _B$ does not depend on the choice of the unital algebra $C$ containing $B,$ so one may
adjoin a unity element to $B$ by one of the most common ways: $C= K \times B,$ $(k,b) (k',b'):=
(k k', kb' + k'b + bb').$
Let  $\phi :A\to B$  be a morphism of $G$-graded algebras, and   $C'= K \times B$ the unital algebra obtained
same way as $C.$  Then   $\phi $ obviously extends to a morphism of unital algebras $C' \to C.$ The homomorphism
$\phi ^\#:A\# G\to B\# G$ between the smash products induces a  $G$-graded
homomorphism $ \tilde {\phi }:  (  A\# G )  \rtimes _{\beta ^A}G\to ( B\#
G)  \rtimes _{\beta ^B}G$, determined by
\[\tilde {\phi }(c^A\delta _t)=\phi ^\#(c^A)\delta _t
=\sum _{r,s\in G}\phi ^\#(c^A_{r^{-1}s}\munit _{r,s}) \delta _t
=\sum _{r,s\in G}\phi (c^A_{r^{-1}s})\munit _{r,s} \delta _t  .\] The map
$\phi \mapsto \tilde {\phi }$ is the morphism level of the functor
$A\mapsto A\# G\rtimes _{\beta ^A}G$. Similarly,
$\phi ^{\textrm {fin}}:{\rm FMat}_G(A)\to {\rm FMat}_G(B) $
% M_G{^\textrm {fin}}(B)$,
given
by $\phi ^{\textrm {fin}}(d)(r,s)=\phi (d(r,s))$, is the morphism level of
the functor $A\mapsto {\rm FMat}_G(A)$.  Clearly, $\phi ^{\textrm {fin}}(d)$ is also a $G$-graded map.  A direct computation shows
that the diagram below commutes:
\[\xymatrix { (A\# G)\rtimes _{\beta ^A}G
\ar@ {->}[r]^-{\psi _A}_-{\cong }\ar [d]_-{\tilde {\phi }}
&{\rm FMat}_G(A)\ar [d]^{\phi ^{\textrm {fin}}}\\
(B\# G)\rtimes _{\beta ^B}G\ar@ {->}[r]_-{\psi _B}^-{\cong }
&{\rm FMat}_G(B)} \]
Now the proof is complete.
\end {proof}
\subsection {Multipliers}
\par

  %We
%remind the concept  of the multiplier algebra $\mult (A)$ of an algebra $A.$
%We use the right hand side notation   for homomorphisms of left
%$A$-modules, while for homomorphisms of right modules  the  usual
%notation shall be used. In particular,  we write $a \mapsto
% a R $ and   $a \mapsto L a$  for
%$R : {}_{A} A  \to {}_{A} A,$  $L : {A}_{A} \to {A}_{A}$ with $ a \in A.$
We recall that the multiplier algebra $\mult (A)$ of an algebra $A$
is the
 set $$\mult (A)= \{(L,R) \in {\rm End}(A_{A}) \times {\rm End}(_{A} A) :
R(a)b = a L(b) \, \mbox {for all}\, a,b
\in A \}$$ with component-wise addition, and multiplication given by
\[(L,R)(L',R'):=(LL',R'R),\quad \forall (L,R),\,(L',R')\in \mult (A).\]
See, for example,  \cite {dokex} or \cite {fd} for details.
For a multiplier $w = (L,R) \in \mult ( A)$ and $a \in A$ we  set $a w = R(a )$ and
$w a = L (a),$ so that
 one always has $(a w ) b = a (w b)$ $(a,b \in A).$    The first
(resp. second) components of the elements of $\mult ( A)$ are called left (resp.
 right) multipliers of $ A$.

 Consider a graded algebra $B=\oplus _{t\in G}B_t$ over the
group~$G$. We denote  by $\mu :B\to \mult (B)$ the natural map, that is, $\mu (b)=(L_b,R_b)$, where
$L_{b}$ and $R_{b}$ are respectively the maps of left and right
multiplication by $b$.
\par A multiplier $w=(L,R)$ of $B$ is said to have degree $t\in G$ if
$w B_s\subseteq B_{ts}$ and $B_s w\subseteq B_{st}$, $\forall s\in
G$. For instance, the multiplier $\mu (b_t)$ defined by $b_t\in
B_t$, is a multiplier of $B$ of degree $t$.
\par Let
$\mult _t(B):=\{ w \in \mult (B): w \textrm { is of degree }t\}$. It
is not hard to see that $\mult _s(B)\mult _t(B)\subseteq \mult _{st}(B)$, $\forall
s,t\in G$, from which it easily follows that $\mult _1(B)$ is a unital algebra
and each $\mult _t(B)$ is a bimodule over $\mult _1(B)$. On the other
hand, since the family $\{B_t\}_{t\in G}$ is linearly independent, it
follows that  the family $\{\mult _t(B)\}_{t\in G}$ is linearly
independent as well. Thus we get a graded algebra
$\mult _{\mathfrak {g}}(B)=\oplus _{t\in G}\mult _t(B)$, which will be
called the \textit {graded multiplier algebra }of $B$. Note that the
natural map $\mu : B\to \mult _{\mathfrak {g}}(B)$ is now a homomorphism
of graded algebras, and $\mu (B)$ is a graded ideal in
$\mult _{\mathfrak {g}}(B)$.

\section {graded-equivalence}\label {sec:GrEquiv}
Suppose $A$ is an associative idempotent algebra.
Consider, in the
category of all right $A$-modules, the full subcategory
mod-$A$ of the unital and torsion-free modules. Thus a right module $M$
over $A$ is in mod-$A$ if and only if $MA=M$ and $mA=0$ implies $m=0$.
This is a Grothendieck category, that is, an abelian category which is
cocomplete and such that direct limits are exact and there exist
generators. In \cite {garsim} the authors characterized the equivalence
of the categories mod-$A$ and mod-$B$ for idempotent algebras $A$ and $B$
in terms of Morita-type theorems: they proved that these
categories are equivalent if and only if there exists a Morita context
$(A,B, {}_A X _B , {}_B Y_A,\tau _A,\tau _B)$, where the modules ${}_AX$, $X_B$, ${}_BY$, $Y_A$ are
unital  and the bimodule maps $\tau _A:{}_AX\otimes _BY_A\to A$ and
$\tau _B:{}_BY\otimes _AX_B\to B$ are surjective (Proposition~2.6 and
Theorem~2.7 of \cite {garsim}).

\par Since we are working with graded algebras, we are interested in
graded Morita contexts. By a graded Morita context between two
idempotent $G$-graded algebras $A=\oplus _{t\in G}A_t$ and
$B=\oplus _{t\in G}B_t$ we mean a Morita context
$(A,B,{}_A X _B , {}_B Y_A, \tau _A,\tau _B),$  where  $X=\oplus _{t\in
G}X_t$ and $Y=\oplus _{t\in G}Y_t$ are
 $G$-graded bimodules,  and  $\tau _{A}(X_r\otimes _B Y_s)\subseteq A_{rs}$
and $\tau _{B}(Y_r\otimes _A X_s)\subseteq B_{rs},$ $\forall r,s,t\in
G$.  Notice that this extends the concept  of a graded Morita context given for the case
of unital rings in \cite {Boisen}. We say that a bimoldule ${}_A X _B$ is unital if both ${}_AX$, $X_B$ are
unital modules. Equivalently, $AXB=X.$
\begin {defn}\label {def:grEquiv}
Let $A=\oplus _{t\in G}A_t$ and $B=\oplus _{t\in G}B_t$ be two
idempotent $G$-graded algebras. We say that they are {\it graded-equivalent} if there exists a graded Morita context
$(A,B,X,Y,\tau _A,\tau _B)$ with unital bimodules $ {}_A X _B , {}_B Y_A$ for
which $\tau _A$ and $\tau _B$ are
surjective.
\end {defn}
%\marginpar {\fac {Does this give an equivalence between categories of graded modules?}}\marginpar {\fac {Is this an equivalence relation?}}

\par It follows from \cite[Theorem 2.6]{Haefner95} that for graded rings with graded local units the concept of a  graded equivalence in  Definition~\ref{def:grEquiv} is tantamount to that considered earlier.

 \par In general we will work with Morita contexts that are contained
  in a graded algebra, that is, $A$, $B$, $X$ and $Y$ will be
  contained (as graded objects) in a
  certain graded algebra $C$, and all the algebraic operations of the context
  will be inherited from the algebra structure of $C$
  (Proposition~\ref {prop:moridalgebras} below shows that we do not
  lose generality in doing so). In particular $\tau _A$ and $\tau _B$
  will be determined by the product of $C$, and we will omit to
  mention them. We will refer to $\mathsf {M}:=(A,B,X,Y)$ as a Morita context in $C$.

\begin {prop}\label {prop:idga}
Let $B=\oplus _{t\in G}B_t$ be a graded algebra. Then $B$ is idempotent
if and only if $B_r=\sum _{s\in G}B_sB_{s^{-1}r}$, $\forall r\in G$.
\end {prop}
\begin {proof}
Just note that
\[
B^2%=(\oplus _{s\in G}B_s)(\oplus _{t\in G}B_t)
=\sum _{s,t\in G}B_sB_t
=\sum _{r\in G}\sum _{st=r}B_sB_t=\bigoplus _{r\in G}\sum _{s\in
  G}B_sB_{s^{-1}r}.\]
\end {proof}

\begin {thm}\label {thm:geq}
    If $B=\oplus _{t\in G} B_t $ is an idempotent $G$-graded algebra, then
    $B$ and $(B\# G)\rtimes _{\color {black} {\beta }^B} G$
are graded-equivalent.
\end {thm}
\begin {proof}
  Given $r,t\in G$, consider the following subsets of $B\#G$: $X_t(r)=B_{r^{-1}t}\munit _{r,t}$,
$X_t=\oplus _{r\in G}X_t(r)$, $Y_t(r)=B_r\munit _{1,r}$, and $Y_t=\oplus _{r\in
  G}Y_t(r)$ (so $Y_t(r)$ and $Y_t$ do not really depend on~$t$). Define now the following subsets of
$(B\# G)\rtimes _{\color {black} {\beta }^B} G$: $X:=\oplus _{t\in G}X_t\delta _t$, $Y=\oplus _{t\in
  G}Y_t\delta _t$, and $B':=\oplus _{t\in
  G}B_t\munit _{1,t}\delta _t$. Note that the map $B'\to B$ given by
$b_t\munit _{1,t}\delta _t\mapsto b_t$ is an isomorphism of graded
algebras, because
 $$(b_s\munit _{1,s}\delta _s)  (b_t\munit _{1,t}\delta _t) =  b_s\munit _{1,s}   \beta _s (b_t\munit _{1,t}) \delta _{st} =
b_s b_t\munit _{1,st}  \delta _{st}\mapsto b_s b_t.$$
 So it is enough to prove that $\mathsf {M}=((B\# G)\rtimes _{\color {black}{\beta }^B}
G,B',X,Y)$ is a graded-equivalence, which implies that
$(B\#G)\rtimes _{\color {black}{\beta }^B} G$ and $B$ are graded-equivalent. Since $X,Y,$ and
$B'$ are graded according to the grading of the crossed product $(B\# G)\rtimes _{\color {black}{\beta }^B} G$, to
see that $\mathsf {M}$ is a graded-equivalence is enough to show that
it is a Morita equivalence.

\par
 We show first that $XY=(B\#G)\rtimes _{\color {black}{\beta }^B} G$. For all $u,r,s,t$ in $G$ we have:
\begin {gather*}
  (X_u(r)\delta _u ) (Y_{u^{-1}t}(s)\delta _{u^{-1}t})
  =X_u(r)\beta _u(Y_{u^{-1}t}(s))\delta _t
  =(B_{r^{-1}u}\munit _{r,u})(B_s\munit _{u,us})\delta _t\\
  =B_{r^{-1}u}B_s\munit _{r,us}\delta _t
  =B_{r^{-1}u}B_{u^{-1}s'}\munit _{r,s'}\delta _t\subseteq (B\#
  G)\delta _t,
\end {gather*}
where $s'=us$. Therefore:
$$(XY)_t = \sum _{u\in G} (X_u \delta _u) (Y_{u\m t} \delta _{u\m t}) =
\sum _{u,r,s} (X_u(r) \delta _u ) ( Y_{u\m t} (s) \delta _{u\m t}).$$ Hence
\begin {gather*}
  XY=\oplus _{t\in G}(XY)_t\delta _t=
\oplus _{t\in G}\sum _{u\in G}\oplus _{r,s\in G}B_{r^{-1}u}B_{u^{-1}s}\munit _{r,s}\delta _t\\
  =\oplus _{t\in G}\oplus _{r,s\in G}\sum _{u\in G}B_{r^{-1}u}B_{u^{-1}s}\munit _{r,s}\delta _t
  =\oplus _{t\in G}(B\# G)\delta _t=(B\# G)\rtimes _\beta G,
\end {gather*} in view of Proposition~\ref {prop:idga}.
We now show that $((B\# G)\rtimes _{\color {black}{\beta }^B} G)X=X$. Given $u,v,r,s,t\in G$:
\begin {gather*}
(B_{r^{-1}s}\munit _{r,s}\delta _u) (X_{u^{-1}t}(v)\delta _{u^{-1}t})
=(B_{r^{-1}s}\munit _{r,s})\beta _u(B_{v^{-1}u^{-1}t}\munit _{v,u^{-1}t})\delta _t\\
=(B_{r^{-1}s}B_{v^{-1}u^{-1}t}\munit _{r,s}\munit _{uv,t})\delta _t
=[s=uv]B_{r^{-1}s}B_{s^{-1}t}\munit _{r,t}\delta _t.
\end {gather*}
Thus, using again Proposition~\ref {prop:idga}, we see that
\begin {gather*}((B\# G)\rtimes _{\color {black}{\beta }^B} G)X=\oplus _{t\in
  G}(\sum _u\sum _{r,s,v} [s=uv] B_{r^{-1}s}B_{s^{-1}t}\munit _{r,t})\delta _t\\=\oplus _{t\in
G}(\sum _rB_{r^{-1}t}\munit _{r,t})\delta _t=\oplus _tX_t\delta _t=X.
\end {gather*}
Let us show that $YX=B'$. Let $u,r,t,s\in G$. Then:
\begin {gather*}
  (Y_u(r)\delta _u) (X_{u^{-1}t}(s)\delta _{u^{-1}t})
=(B_r\munit _{1,r})\beta _u(B_{s^{-1}u^{-1}t}\munit _{s,u^{-1}t})\delta _t\\
=(B_r\munit _{1,r})(B_{s^{-1}u^{-1}t}\munit _{us,t})\delta _t
=[r=us]B_rB_{r^{-1}t}\munit _{1,t}\delta _t.
\end {gather*}
Hence  $YX=\oplus _t(\sum _{u}(\sum _{r,s} [r=us] B_rB_{r^{-1}t}\munit _{1,t}))\delta _t
=\oplus _tB_t\munit _{1,t}\delta _t=B'$.
\par Now it is easy to check that $XB'=X$: \[XB'=X(YX)=(XY)X=((B\# G)\rtimes _\beta G) X=X. \]
  For the equality $B' Y = Y$ notice first that
\begin {gather*}
 (B_u \munit _{1,u} \delta _u ) (  Y_{u^{-1} t}(r) \delta _{u^{-1}t} ) =
(B_u \munit _{1,u} \delta _u ) (   B_r  \munit _{1,r}  \delta _{u^{-1}t} )\\
=(B_u \munit _{1,u}  ) \beta _u (   B_r  \munit _{1,r}  \delta _{t} )
=B_u \munit _{1,u}      B_r  \munit _{u,ur}  \delta _{t} =
B_u      B_r  \munit _{1,ur}  \delta _{t}.
\end {gather*} Then $ B' Y = \oplus _t ( \sum _u \sum _r B_u B_r  \munit _{1,ur} )  \delta _{t}
=  \oplus _t ( \sum _{s, u}  B_u B_{u^{-1}s}  \munit _{1,s} )  \delta _{t} ,$ where $ s=ur.$ With one more use of  Proposition~\ref {prop:idga} this shows that  $B' Y = Y.$ Finally,
\[Y ( (B\# G)\rtimes _{\color {black}{\beta }^B} G) = Y (X Y) =  (Y X)  Y = B' Y =Y,  \]
which ends the proof.  \end {proof}

\par If $C$ is an algebra and $e\in \mult (C)$ is an idempotent element, then
$C':=e C e$ is clearly an algebra.
  Suppose in addition that $C=\oplus _{t\in
  G}C_t$ is graded over $G$, and $e$ is a multiplier of $C$ of degree
1. Define $C'_t:=eC_te$. Then $C'_sC'_t=eC_seC_te\subseteq
eC_sC_te\subseteq eC_{st}e=C'_{st}$, so it follows that $C'$ is also a
graded algebra, because $C'=\oplus _{t\in G}C'_t$.

\begin {prop}\label {prop:moridalgebras}
Let $A=\oplus _{t\in G}A_t$ and $B=\oplus _{t\in G}B_t$ be idempotent
algebras graded over the group $G$. Then the following statements are
equivalent:
\begin {enumerate}
\item The algebras $A$ and $B$ are graded-equivalent.
\item There exist an idempotent graded algebra $C=\oplus _{t\in G}C_t$
  and $e=e^2\in \mult _1(C)$ (the algebra of multipliers of degree 1 of
  $C$), such that: $A\cong eCe$ and $B\cong (1-e)C(1-e)$ as graded
  algebras, and $CeC=C=C(1-e)C$.
\end {enumerate}
\end {prop}
\begin {proof}
Suppose first that $(A,B,X,Y,\tau _A,\tau _B)$ is a graded Morita
context between $A$ and $B$. Let
$\mathbb {L}:=\{\begin {pmatrix}a&x\\y &b\end {pmatrix}: a\in A,b\in
B, x\in X, y\in Y\}$ with entry-wise addition, and the product given by
the Morita context, as follows:
\[\begin {pmatrix}a&x\\y &b\end {pmatrix}
\begin {pmatrix}a'&x'\\y '&b'\end {pmatrix}
=\begin {pmatrix}aa'+\tau _A(x\otimes
  y)&ax'+xb'\\ya '+by'&\tau _B(y\otimes x')+bb'\end {pmatrix}.
\]
These operations give an associative algebra structure on $\mathbb {L}$,
as it is easy to check. Moreover, since the algebras $A$ and $B$ are
idempotent, and the modules of the Morita context are unital, it
follows that $\mathbb {L}$ is
also an idempotent algebra. If
\[\mathbb {L}_t:=\{\begin {pmatrix}a_t&x_t\\y _t&b_t\end {pmatrix}: a_t\in
A_t,b_t\in B_t, x_t\in X_t, y_t\in Y_t\},\] then the fact that the
Morita context is graded implies $\mathbb {L}_s\mathbb {L}_t\subseteq
\mathbb {L}_{st}$, $\forall s,t\in G$. On the other hand it is clear that
$\mathbb {L}=\oplus _{t\in G}\mathbb {L}_t$, so $\mathbb {L}$ is a graded
algebra over $G$.
\par Now let $L,R:\mathbb {L}\to \mathbb {L}$ be the
maps given by
$L(\begin {pmatrix}a&x\\ y&b\end {pmatrix})
=\begin {pmatrix}a&x\\0&0\end {pmatrix}$ and
$R(\begin {pmatrix}a&x\\y &b\end {pmatrix})
=\begin {pmatrix}a&0\\y &0\end {pmatrix}$. Then $L(\mathbb {L}_t)\subseteq
\mathbb {L}_t$, $R(\mathbb {L}_t)\subseteq \mathbb {L}_t$, and routine
matrix computations
show that $L(cc')=L(c)c'$, $R(cc')=cR(c')$, and $cL(c')=R(c)c'$,
$\forall c,c'\in \mathbb {L}$. That is, $e:=(L,R)$ is a multiplier of
$\mathbb {L}$, of degree 1. It is clear that $e^2=e$. Note
that we may conveniently think of $e$ as the
matrix $\begin {pmatrix}1&0\\0&0\end {pmatrix}$ acting in the obvious
way on $\mathbb {L}$: $L$ is multiplication on the left by this matrix, while
$R$ corresponds to multiplication on the right. Now,
$e\mathbb {L}_te=R(L(\mathbb {L}_t))=
\begin {pmatrix}A_t&0\\0&0\end {pmatrix}$, and
$(1-e)\mathbb {L}_t(1-e)=(1-R)((1-L)(\mathbb {L}_t))
=\begin {pmatrix}0&0\\0&B_t\end {pmatrix}$,
which are obviously isomorphic to $A_t$ and $B_t$ respectively. On the
other hand:
\[\mathbb {L}e\mathbb {L}=\begin {pmatrix}A&X\\Y &B\end {pmatrix}
\begin {pmatrix}A&X\\0&0\end {pmatrix}
=\begin {pmatrix}A^2&AX\\YA &\tau _B({}_BY\otimes _AX_B)\end {pmatrix}
=\mathbb {L},\]
where the latter equality is due to the facts that $A$ is idempotent,
${}_AX$ and $Y_A$ are unital, and $\tau _B$ is surjective. In a similar
way we conclude that $\mathbb {L}(1-e)\mathbb {L}=\mathbb {L}$.
\par
Conversely, suppose $C=C^2$ and $e\in \mult _1(C)$ are such that
$A':=eCe\cong A$ and $B':=(1-e)C(1-e)\cong B$ as graded algebras, and
$CeC=C=C(1-e)C$. Let
$X:=eC(1-e)$, $X_t:=eC_t(1-e)$, $Y:=(1-e)Ce$, and $Y_t:=(1-e)C_te$. It
is clear that $X$ is an
$(A',B')$-bimodule. Moreover:
\begin {gather*}
A'_sX_t=(eC_se)(eC_t(1-e))\subseteq
eC_sC_t(1-e)\subseteq X_{st},  \\
A'X=(eCe)(eC(1-e))=e(Ce^2C)(1-e)=eC(1-e)=X.
\end {gather*}
Similarly we have $X_sB'_t\subseteq X_{st}$ and $XB'=X$, so
${}_AX$ and $X_B$ are unital and graded modules. In the same way we
conclude that ${}_BY$ and $Y_A$ are
unital and graded modules. Since
$$XY=(eC(1-e))((1-e)Ce)=e(C(1-e)C)e=eCe=A',$$ and $YX=B'$
after a similar computation, we have that the natural maps
$\tau _{A'}:X\otimes _{B'}Y\to A'$ and $\tau _{B'}:Y\otimes _{A'}X\to B'$
associated to the multiplication on $C$ are both surjective, and it is
easily seen that $\tau _{A'}(X_s\otimes _{B'}Y_t)\subseteq A'_{st}$ and
$\tau _{B'}(Y_s\otimes _{A'}X_t)\subseteq B'_{st}$, $\forall s,t\in
G$. Then $(A',B',X,Y,\tau _{A'},\tau _{B'})$ is a graded-equivalence between
$A'$ and $B'$. Since $A\cong A'$ and
$B\cong B'$ as graded algebras, we are done.
\end {proof}

\begin {defn}
Let $\mathsf {M}=(A,B,X,Y,\tau _A,\tau _B)$ be a graded Morita context giving a
graded-equivalence between the idempotent graded algebras $A$ and
$B$. The graded algebra $\mathbb {L}$ constructed out of this Morita
context as in the first part of the proof of
Proposition~\ref {prop:moridalgebras}
will be called the {\it graded Morita algebra} (or the {\it graded linking algebra}) of  $\mathsf {M}$.
We will write $\mathbb {L}(\mathsf {M})$ if we need to stress the dependence of $\mathbb {L}$ on the
Morita context $\mathsf {M}.$
%, and we will also refer to $\mathbb {L}(\mathsf {M})$ as the {\it Morita algebra} of the context $\mathsf {M}$.
\end {defn}
\par We end the section by showing that graded-equi\-va\-len\-ce is an
  equi\-va\-len\-ce relation:
\begin {prop}\label {prop:ge is eq}
  Graded-equivalence is an equivalence relation
for graded idempotent algebras.
\end {prop}
\begin {proof}
Suppose $\mathsf {M}=(A,A',X,Y,\tau _A,\tau _{A'})$ and
$\mathsf {N}=(A',B,X',Y',\tau '_{A'},\tau _B)$ are graded Morita contexts
giving graded-equivalences between $A$ and $A'$ and between $A'$ and
$B$ respectively.
%Without loss of generality we suppose that the sets
%of the context $\mathsf {M}$ are subsets of an algebra $C$, and that
%$\tau _A$ and $\tau _{A'}$ are given by the
%multiplication of $C$. A similar assumption we make about the context
%$\mathsf {N}$.
 As in \cite [p. 30]{nastoyst} we consider $_A\bar {X}_B = X\otimes _{A'} X'$ as a
$G$-graded bimodule, whose $t$-homogeneous component     $\bar {X}_t, $ $(t\in G)$ is
 the $K$-submodule  of $\bar {X}$ generated by all elements $x \otimes x',$
$x \in X_r, x'\in X'_s,$ such that $rs=t.$
Thus, for  $t\in G,$      $\bar {X}_t = \sum _{s\in G}{X_s}\otimes _{A'} X'_{s\m t}.$
Similarly, $_B\bar {Y}_A = Y'\otimes _{A'}Y$ is a $G$-graded bimodule, with
 $\bar {Y}_t  = \sum _{s\in G}{Y'_s}\otimes _{A'} Y_{s\m t}.$ Obviously,   $_A\bar {X}_B$  and  $_B\bar {Y}_A $ are
unital bimodules.
 %Note first that the $(A,B)$-bimodule $\bar {X}:=X\otimes _{A'_1} X'$ is
%$G$-graded. In fact, if $X=\oplus _{r\in G}X_r$ and $X'=\oplus _{s\in
  %G}X'_s$, we have $X\otimes _{A'_1} X'=\oplus _{r,s\in
 % G}X_r\otimes _{A'_1}X'_s=\oplus _{t\in G}\bar {X}_t$, where
%$\bar {X}_t:=\oplus _{rs=t}X_r\otimes _{A'_1}X'_s$. It is clear that,
%with the grading $\bar {X}=\oplus _{t\in G}\bar {X}_t$, $\bar {X}$ becomes
%a unital graded $(A,B)$-bimodule. Similarly, we endow
%$\bar {Y}:=Y'\otimes _{A'_1}Y$ with a structure of graded
%$(B,A)$-bimodule.
Denote by $\rho _A$ the following composition of
surjective $(A,A)$-bimodule maps:
\begin {gather*}
\bar {X}\otimes _B\bar {Y}
=(X\otimes _{A' } X')\otimes _B(Y'\otimes _{A' }Y)\\
\cong X\otimes _{A' } (X'\otimes _BY')\otimes _{A' }Y
\stackrel {id\otimes \tau '_{A'}\otimes id}{\longrightarrow }
X\otimes _{A' } A'\otimes _{A' }Y  \\
%\stackrel {id\otimes id\otimes id}{\longrightarrow }X\otimes _{A'}
%A'\otimes _{A'}Y
 \longrightarrow X\otimes _{A'}Y
\stackrel {\tau _{A}}{\longrightarrow } A.
\end {gather*}
Note that $\rho _A(\bar {X}_s\otimes _B\bar {Y}_t)\subseteq A_{st}$,
$\forall s,t\in G$. In fact, if $u,v\in G$ and
$x\in X_u, x'\in X'_{u\m s}, y' \in Y'_{v}, y\in Y_{v\m t,}$ we have:
\[\rho _A((x\otimes
x')\otimes (y'\otimes y))=
\tau _A ( (x\otimes {\tau }'_{A'}( x'\otimes y')y) \in \tau _A(X_u\otimes Y_{u^{-1}st})\subseteq
A_{st}.\]
In the same way we construct a surjective $(B,B)$-bimodule map
$$\rho _{B}:\bar {Y}\otimes _{A}\bar {X}\to B,$$  such that
$\rho _B(\bar {Y}_s\otimes _{\color {black} A}\bar {X}_t)\subseteq B_{st}$, $\forall s,t\in
G$.  It then follows that $\mathsf {MN}:=(A,B,\bar {X},\bar {Y},\rho _A,\rho _B)$ is a
 graded-equivalence between $A$ and $B$.
\end {proof}

\begin {rem}\label {rem:BFMat}
Note that combining Proposition~\ref {prop:ge is eq} with
Theorem~\ref {thm:duality} and Theorem~\ref {thm:geq}, we obtain that
any $G$-graded idempotent algebra $B$ is graded-equivalent to ${\rm
  FMat}_G(B)$. \end {rem}

\section {Strong-graded-equivalence}\label {sec:StrGrEquiv}
\par If $A=\oplus _{t\in G}A_t$ is a $G$-graded algebra, then for each
$t\in G$ the set $D_t^A:=A_tA_{t^{-1}}$ is a two-sided ideal of
$A_1$.
\par Suppose that $(A,B,X,Y,\tau _A,\tau _B)$ is a graded Morita context
between idempotent $G$-graded algebras $A$ and $B$. Since $D_{t^{-1}}^B$ is a
subalgebra of $B$, we have a natural map
 $\mu _A^t:X_t\otimes _{D_{t^{-1}}^B}Y_{t^{-1}}\to
X\otimes _{B}Y$ (observe that whenever $X_t$ or $Y_{t^{-1}}$ are unital
  $D_{t^{-1}}^B$-modules, then
  $Y_t\otimes _{D_{t^{-1}}^B}X_{t^{-1}}=Y_t\otimes _{B_1}X_{t^{-1}}$).
 Thus composing this map with $\tau _A$ we
obtain an $A_1$-bimodule map
$\tau _A^t:=\tau _A\mu _A^t:X_t\otimes _{D_{t^{-1}}^B}Y_{t^{-1}}\to A_1$.
Similarly, we have a natural map
$\mu _B^t:Y_t\otimes _{D_{t^{-1}}^A}X_{t^{-1}}\to
Y\otimes _{A}X$, and also a $B_1$-bimodule  map
$\tau _B^t:=\tau _B\mu _B^t:Y_t\otimes _{D_{t^{-1}}^A}X_{t^{-1}}\to B_1$.
\par Suppose that $D_t^A$ and $D_t^B$ are idempotent, $\forall t\in
G$. We will say that the Morita context above is $\textit {strong}$, if
each $X_t$ is a unital $(D_t^A, D_{t^{-1}}^B)$-bimodule and  each $Y_t$ is a unital
$(D_t^B, D_{t^{-1}}^A)$-bimodule.  Under this condition observe that, for each $t\in G$, the ranges of
$\tau _A^t$ and $\tau _B^t$ are contained in $D_t^A$ and $D_t^B$
respectively.  Indeed,
 $$
  \tau _A^t(X_t\otimes _{D_{t^{-1}}^B}Y_{t^{-1}}) =
  \tau _A^t(D_t^A X_t\otimes _{D_{t^{-1}}^B}Y_{t^{-1}}) =$$ $$=
  D_t^A \tau _A^t(X_t\otimes _{D_{t^{-1}}^B}Y_{t^{-1}})   \subseteq
  D_t^A A_1 \subseteq D_t^A,
  $$
  %$$
  %\tau _A^t(X_t\otimes _{D_{t^{-1}}^B}Y_{t^{-1}}) =
  %\tau _A^t(D_t^AX_t\otimes _{D_{t^{-1}}^B}Y_{t^{-1}}) =$$ $$=
  %\tau _A^t(A_tA_{t^{-1}}X_t\otimes _{D_{t^{-1}}^B}Y_{t^{-1}}) =
  %A_t\tau _A^t(A_{t^{-1}}X_t\otimes _{D_{t^{-1}}^B}Y_{t^{-1}}) \subseteq $$ $$ \subseteq
  %A_t\tau _A^t(X_1\otimes _{D_{t^{-1}}^B}Y_{t^{-1}}) \subseteq
  %A_tA_{t^{-1}} = D_t^A,
  %$$
  and similarly for the range of $\tau _B^t$.

The following definition is an algebraic adaptation of \cite [Definition~2.6]{abf}, which in turn is a generalization
of the notion of Morita-Rieffel equivalence of partial actions introduced in \cite {abenv}. As we shall see later in Theorem~\ref {thm:moritaskew}, Definition~\ref {defn:sgeq} is a generalization of Morita equivalence of
\textit {product partial actions}, which are defined in Section~\ref {sec:ppa}.

\begin {defn}\label {defn:sgeq}
Let $A=\oplus _{t\in G}A_t$ and $B=\oplus _{t\in G}B_t$ be two idempotent  graded
algebras such that $D_t^A$ and $D_t^B$ are idempotent, $\forall t\in
G$. We say that $A$ and $B$ are \textit {strongly-graded-equivalent} if there exists a strong-graded Morita context
$(A,B,X,Y,\tau _A,\tau _B)$ with surjective $\tau _A$ and $\tau _B,$  such that $\tau _A^t(X_t\otimes _{D_{t^{-1}}^B}
Y_{t^{-1}})=D_t^A$ and $\tau _B^t(Y_t\otimes _{D_{t^{-1}}^A} X_{t^{-1}})=D_t^B$,
$\forall t\in G$.
\end {defn}

\par Evidently,  if $(A,B,X,Y,\tau _A,\tau _B)$
is a  strong-graded Morita context, then the bimodules $ {}_A X _B , {}_B Y_A$ are unital. Hence, strongly-graded-equivalent algebras are graded-equivalent.   Moreover,  notice that
if $D_t^A$ is idempotent then $D_t^A$ is a unital $(A_1,
A_1)$-bimodule, since, for example,    $D_t^A = (D_t^A)^2 \subseteq
D_t^AA_1 \subseteq D_t^A.$

\begin {prop}\label {prop:strong}
  Let $\mathsf {M}=(A,B,X,Y,\tau _A,\tau _B)$ be a strong-graded Morita context.
Then:
  \begin {enumerate}
  \item $A_rX_1=X_r=X_1B_r$ and $B_rY_1=Y_r=Y_1A_r$, $\forall r\in G$.
  \item $D_r^AX_1=X_1D_r^B$ and $D_r^BY_1=Y_1D_r^A$, $\forall r\in G$.
  \item $\tau _A(X_r\otimes _B  Y_s)=A_r\tau _A(X_1\otimes _B Y_1)A_s\subseteq
    A_rA_s$ and $\tau _B(Y_r\otimes _AX_s)=B_r\tau _B(Y_1\otimes _A X_1)B_s\subseteq B_rB_s$, $\forall r,s\in G$.
  \item If the context $\mathsf {M}$ is a strong-graded-equivalence and
    each $A_t$ is a unital $A_1$-bimodule and each $B_t$ is a unital
    $B_1$-bimodule,  then $\tau _A(X_r\otimes _B  Y_s)=A_rA_s$ and
    $\tau _B(Y_r\otimes _AX_s)=B_rB_s$.
  \end {enumerate}
 \end {prop}
 \begin {proof} Since $X_r$ is a unital left $D_r^A$-module, we
    have \[A_rX_1\subseteq X_r=D_r^AX_r =A_rA_{r^{-1}}X_r \subseteq
    A_rX_1,\] so $A_rX_1=X_r$. The other
    identities in (1) are proved similarly, and (2) follows at once
    from (1). Now, since $\tau _A$ is an $A$-bimodule map:
\[\tau _A(X_r\otimes Y_s)=\tau _A(A_rX_1\otimes Y_1A_s)
=A_r\tau _A(X_1\otimes Y_1)A_s\subseteq A_rA_1A_s \subseteq A_rA_s.\]
An analogous argument proves the second claim of (3). Finally, (4) follows from (3),
because %\todo {Here $A_rA_1A_s=A_rA_s$ is used}
the assumption implies $\tau _A(X_1\otimes _B Y_1)=A_1$ and $\tau _B(Y_1\otimes _A X_1)=B_1$.
    \end {proof}

 \par Note that, if in Proposition~\ref {prop:strong} the context
 $(A_1,B_1,X_1,Y_1,\tau _A^1,\tau _B^1)$ is a Morita equivalence, and each $A_t$ is a unital
$A_1$-bimodule and each $B_t$ is a unital $B_1$-bimodule,  then
 both of the inclusions in (3) above are actually equalities. So we
 have:

\begin {cor}\label {cor:sgt-sgeq}  Let $A = \oplus _{t\in G } A_t$ and $B= \oplus _{t\in G} B_t$ be
idempotent $G$-graded algebras such that  each   $A_t$ is a
unital $A_1$-bimodule and each $B_t$ is a unital $B_1$-bimodule.
  Then a  strongly-graded Morita context
  $(A,B,X,Y,\tau _A,\tau _B)$ is a strong-graded-equivalence between $A$
  and $B$ if and only if the Morita context
 $(A_1,B_1,X_1,Y_1,\tau _A^1,\tau _B^1)$ is a Morita equivalence.
\end {cor}
\begin {proof} The `only if' part  is evident and for the `if' part it remains to show that if
$(A_1,B_1,X_1,$ $ Y_1,\tau _A^1,\tau _B^1)$ is a Morita equivalence then
$\tau _A$ and $\tau _B $ are surjective. Indeed, using (1) of Proposition~\ref {prop:strong} we see, for instance,  that $\tau _A$ is surjective as follows:
\begin {gather*}\tau _A (X \otimes _B Y )  = \sum _{r,s \in G} \tau _A (X_r \otimes _B Y_s)  =
\sum _{r,s \in G} A_r \tau ^1 _A (X_1 \otimes _{B_1} Y_1) A_s \\
= \sum _{r,s \in G} (A_r A_1) (A_1 A_s) =A^2= A,
\end {gather*} thanks to the unital condition on each $A_s.$  Analogously, we show that
$\tau _A^t(X_t\otimes _{D_{t^{-1}}^B}
Y_{t^{-1}})=D_t^A$ and $\tau _B^t(Y_t\otimes _{D_{t^{-1}}^A} X_{t^{-1}})=D_t^B$,
$\forall t\in G$.  \end {proof}

\begin {cor}\label {cor:sgt-sgeq2}
If $A$ and $B$ are strongly-graded-equivalent, then the algebras
$D_t^A,D_{t^{-1}}^A, D_t^B$ and $D_{t^{-1}}^B$ are Morita
equivalent to each other.
\end {cor}
\begin {proof}  Note that if
  $(A,B,X,Y,\tau _A,\tau _B)$ gives a strong-graded-equivalence between
  the graded algebras $A$ and $B,$ then  by definition the context
  $(D_t^A,D_{t^{-1}}^B,X_t,Y_{t^{-1}},\tau _A^t,\tau _B^{t^{-1}})$ is a Morita
  equivalence, $\forall t\in G$.
 Since Morita equivalence of idempotent algebras is transitive,
  we only need to show that $D_t^A$ and $D_t^B$ are Morita
equivalent, $\forall t\in G$. But Proposition~\ref {prop:strong} shows
that if $(A,B,X,Y,\tau _A,\tau _B)$ is a strong-graded-equivalence between $A$ and $B$, then $D_t^AX_1$ is a unital
$(D_t^A,D_{t}^B)$-bimodule and $D_t^BY_1$ is a unital
$(D_t^B,D_{t}^A)$-bimodule. Moreover:
\begin {gather*}
\tau _A(D_t^AX_1\otimes _BD_t^BY_1)
=\tau _A(D_t^AX_1D_t^B\otimes _BY_1)\\
=D_t^A\tau _A(X_1\otimes _BY_1)
=D_t^A.\end {gather*}
Similar computations show that also
$\tau _B(D_t^BY_1\otimes _AD_t^AX_1)=D_t^B$. Thus
$(D_t^A,D_{t}^B, D_t^AX_1$ $,D_t^BY_1,\tilde {\tau }_A,\tilde {\tau }_B)$
defines a Morita equivalence between $D_t^A$ and $D_{t}^B,$ where $\tilde {\tau }_A$ is the composition of
the natural map
	$ D_t^AX_1 \otimes _{ D_t^{B}  } D_t^B  Y_1 \to X \otimes _{B} Y $ with $\tau _A,$ and $\tilde {\tau }_B$ is
defined symmetrically.
\end {proof}

\par Observe that every graded algebra $A=\oplus _{t\in G}A_t$ defines
a trivial Morita
context, namely, the context $(A,A,A,A)$, where the range maps are
given by the product of $A$.
 Note that this context is  a strong-graded-equivalence
between $A$ and itself if and only if $A$ satisfies the following
property:
\begin {equation}\label {eq:prp}
    A_r=A_rA_{r^{-1}}A_r,\quad \forall r\in G.
\end {equation} In particular, \eqref {eq:prp}  implies that $A$ is an idempotent algebra, as
$A^2\supseteq A_tA_{t^{-1}}A_t = A_t$ for each $t\in G.$ Moreover, it also follows from  \eqref {eq:prp}
that the algebra $A_1$ is  idempotent, as well as each ideal $D_t^A$ of $A_1,$ $(t\in G).$ In addition,
$A_1 A_t \supseteq A_t A_{t^{-1}} A_t = A_t$ and consequently,  $A_1A_t = A_t, $ $t\in G.$ Similarly,
  $A_tA_1 = A_t, $ so that each  $A_t$ is  a unital  $A_1$-bimodule. As a consequence, we see that
each $D_t^A$ is also a unital   $A_1$-bimodule.

\par Observe that if $A$ satisfies \eqref {eq:prp} then $A$ satisfies
each of the following two conditions
\begin {equation}\label {eq:prp2}
A_{r^{-1}}A_{r}A_s=A_{r^{-1}}A_{rs} \quad \text {and}  \quad A_sA_{r}A_{{r^{-1}}}=
A_{sr}A_{r^{-1}},\quad \forall r,s\in G.
\end {equation} Indeed,
%on the one hand it is clear that the latter condition implies \eqref {eq:prp}. On the other hand,
if $A$ satisfies \eqref {eq:prp}, then:
\[A_{r^{-1}}A_{rs}=A_{r^{-1}}A_{r}A_{r^{-1}}A_{rs}
  \subseteq A_{r^{-1}} A_rA_{r^{-1}rs}
  = A_{r^{-1}} A_rA_{s}
  \subseteq A_{r^{-1}}A_{r s}.
\] The second equality is proved similarly. Observe that if each $A_t$ is a unital left (or right)
$A_1$-module, then  \eqref {eq:prp2} implies  \eqref {eq:prp}.

Notice that equalities of the form \eqref {eq:prp2} appear in the
definition of the concept of a partial representation of a group (see
for example \cite [Definition~2.1]{DEP}). As an effect, if $A$ satisfies \eqref {eq:prp} then,
for all  $t,s \in G$ we have:
\begin {equation}\label {eq:prp3} A_t D_s^A= D_{ts}^A A_t, \;\;\;\;\;  D_t^A D_s^A = D_s^A D_t^A,
\quad \forall t,s\in G,
\end {equation}  which are   the analogues of  useful properties
of partial representations (see (2) and (3) in \cite {DEP}). One can readily adapt the easy computations in
\cite {DEP}, and for an illustration  we  prove  the first equality,
 from which the second one follows easily:
\begin {gather*}
 A_t D_s^A =  A_t A_s A_{s\m } =    A_{ts} A_{s\m } =  A_{ts} A_{s\m
   t\m }   A_{ts}  A_{s\m }=\\
  A_{ts}A_{s\m t\m }A_t  =D_{ts}^AA_t.
\end {gather*}

Note also that  \eqref {eq:prp} is one of the  key properties which characterize
graded algebras as crossed products by twisted partial actions (see \cite [Theorem 6.1]{des}).
 Recall that $A=\oplus _{t\in G}A_t$ is called strongly-graded if
$A_rA_t=A_{rt}$ for all $r,s,\in G.$ Obviously, each strongly-graded algebra satisfies   \eqref {eq:prp},
so we give the next:
\begin {defn}\label {defn:prp}
A graded algebra $A=\oplus _{t\in G}A_t$ over the group $G$ sa\-tis\-fy\-ing
\eqref {eq:prp} is said to be {\it partially-strongly-graded}.
%\marginpar {\fac {Does the strong gr. equivalence define an eq. relation in the class of gr. algs with the prp?}}
\end {defn}
For instance, the so called epsilon-strongly-graded algebras con\-si\-de\-red in \cite {nop} are
partially-strongly-graded.

%\par %It is clear that if $A=\oplus _{t\in G}A_t$ a graded algebra with
%$D_t^A$ is idempotent, and therefore $A$ is  automatically
%idempotent. Consequently,
%Notice that  in Proposition~\ref {prop:strong}, Corollary~\ref {cor:sgt-sgeq} and
%Corollary~\ref {cor:sgt-sgeq2}, the hypotheses on the algebras $A$ and
%$B$ are fulfilled whenever the algebras $A$ and $B$ are supposed {\color {red} to be partially-strongly-graded.}
\par A graded algebra $A$ is strongly-graded precisely when
$D_t^A=A_1$ and $A_1A_t=A_t$
$\forall t\in G$, that is, whenever each $(A_t,A_{t^{-1}})$ defines a
Morita autoequivalence of $A_1$. Similarly, $A$ is  partially-strongly-graded exactly when each
$(D_t^A,D_{t^{-1}}^A,A_t,A_{t^{-1}})$ is a Morita equivalence
between $D_t^A$ and $D_{t^{-1}}^A$ (we omitted the trace maps, which
are given by the product of $A$).

\begin {comment}
\fac {Let $X=\oplus _{t\in G}X_t$ be a graded right module over the
partially-strongly-graded algebra $B=\oplus _{t\in G}B_t$. If $r$ is
such that $X_rD_{r^{-1}}^B=X_r$, then $X_r=X_rD_{r^{-1}}^B\subseteq X_1B_r\subseteq
X_r$, so $X_r=X_1B_r$. Conversely, if $X_r=X_1B_r$, we
have $X_r=X_1B_r=X_1B_rB_{r^{-1}}B_r\subseteq X_rD_{r^{-1}}^B\subseteq
X_r$, which implies $X_r=X_rD_{r^{-1}}^B$. Therefore the fact that each $X_r$
is a unital $D_{r^{-1}}^B$-module can be expressed equivalently by the
condition $X_1B_r=X_r$, $\forall r\in G$. Similarly, if $X$ is a
graded left $B$-module, one sees that $D_r^BX_r=X_r$ if and only if
$B_rX_1=X_r$.
}
\end {comment}
The next result extends formul\ae \ \eqref {eq:prp2} and
  \eqref {eq:prp3} above for graded modules over partially-strongly-graded
  algebras:
\begin {prop}\label {prop:prp} Let $B=\oplus _{t\in G}B_t$ be a partially-strongly-graded algebra, and $X=\oplus _{t\in G}X_t$ a graded module
  over $B$. Then:
  \begin {enumerate}
     \item $X$ is partially-strongly-graded, in the sense that for all
       $r,s\in G$ we have $X_{r}B_sB_{s^{-1}}=X_{rs}B_{s^{-1}}$, if
       $X$ is a right $B$-module, and
       $B_{r^{-1}}B_rX_{s}=B_{r^{-1}}X_{rs}$ if $X$ is a left
       $B$-module
     \item If $X$ is a right module, then $X_r$ is a unital right
       $D_{r^{-1}}^B$-module if and only if $X_r=X_1B_r$, and if $X$
       is a left module, then $X_r$ is a unital left $D_{r}^B$-module
       if and only if $X_r=B_rX_1$.
     \item If $A=\oplus _{t\in G}A_t$ is a partially-strongly-graded
       algebra and $X=\oplus _{t\in G}X_t$ is a graded
       $(A, B)$-bimodule such that each $X_r$ is a unital
       $(D_r^A, D^B_{r^{-1}})$-bimodule, $\forall r\in G$, then
       $X_rD^B_s=D_{rs}^AX_r$, $\forall r\in G$.
  \end {enumerate}
\end {prop}
\begin {proof} Suppose
  that $X$ is a graded right $B$-module. Since $B$
  is strongly-graded we have
  $X_{rs}B_{s^{-1}}=X_{rs}B_{s^{-1}}B_sB_{s^{-1}}\subseteq
  X_rB_sB_{s^{-1}}\subseteq X_{rs}B_{s^{-1}}.$ In particular we get
  $X_rD_{r^{-1}}^B=X_1B_r$, so the first two claims are proved for a
  right $B$-module $X$. A similar argument for
  a left module $X$ concludes the proof of (1) and (2). As for (3)
  note that by (1) we have $X_rD^B_s =X_{rs}B_{s^{-1}}$. On the other
  hand, by (2) we have $A_rX_1=X_r=X_1B_r$, $\forall r\in G$. Thus
$X_{rs}B_{s^{-1}}=A_{rs}X_1B_{s^{-1}}
=A_{rs}A_{s^{-1}}X_1
=D^A_{rs}A_{rs}A_{s^{-1}}X_1
\subseteq D_{rs}^AX_r$. Hence $X_rD^B_s\subseteq D_{rs}^AX_r$. A
symmetric argument shows that
$D_t^AX_r\subseteq X_rD_{r^{-1}t}$, $\forall r,t\in G$. This ends the
proof, since in particular we get $D_{rs}^AX_r\subseteq
X_rD^B_{r^{-1}rs}=X_rD_{s}^B$.
\end {proof}

% FULLY UNITAL

% \hrule

% \begin {defn}\label {defn:fully}
% Let $X$ be a graded right (left) module over the partially-strongly-graded algebra $B$. If $X$ is such that $X_r$ is a unital right
% $D^B_{r^{-1}}$-module (respectively a unital left
% $D^B_{r}$-module) for all $r\in G$, we will say that $X$ is a
% fully unital graded $B$-module. If $X$ is a graded $(A, B)$-bimodule over
% the partially-strongly-graded algebras $A$ and $B$, we will say that
% $X$ is fully unital if it is both fully unital as a left $A$-module
% and as a right $B$-module.
% \end {defn}
% For
% instance, if $(A,B,X,Y,\tau _A,\tau _B)$ is a strong Morita context,
% where $A$ and $B$ are partially-strongly-graded algebras, then $X$ is
% a fully unital $(A, B)$-bimodule, and $Y$ is a fully unital
% $(B, A)$-bimodule.

  % Notice that, by (2) of Proposition~\ref {prop:prp}, a right (left)
% graded module $X$ over the partially-strongly-graded algebra $B$ is fully unital if and only
% if $X_r=X_1B_r$ (respectively: $X_r=B_rX_1$) $\forall r\in G$.
% \begin {rem}\label {rem:fully}
% If $X$ is a graded right module over the partially-strongly-graded algebra $B$, and we define $Z_r:=X_1B_r$,
% $\forall r\in G$, then it is easily checked that $Z:=\oplus _{t\in
% G}Z_t$ is a fully unital right $B$-module.
% \end {rem}

% \hrule

\begin {prop}\label {prop:strongmoridalgebras}
Let $A=\oplus _{t\in G}A_t$ and $B=\oplus _{t\in G}B_t$ be partially-strongly-graded
algebras.  Then the following
statements are equivalent:
\begin {enumerate}
\item The algebras $A$ and $B$ are strongly-graded-equivalent.
\item There exist a partially-strongly-graded algebra $C=\oplus _{t\in G}C_t$  and $e=e^2\in \mult _1(C)$ as in (2)
  of Proposition~\ref {prop:moridalgebras}, which moreover
  satisfy, $\forall t\in G$:
   \begin {enumerate}
     \item $eD_t^Ce\cong D_t^A$, $(1-e)D_t^C(1-e)\cong D_t^B$
     \item $D_t^CeD_t^C=D_t^C=D_t^C(1-e)D_t^C$.
   \end {enumerate}
\end {enumerate}
\end {prop}
\begin {proof}
Let $\mathbb {L}=\oplus _{t\in G}\mathbb {L}_t$ be the Morita algebra of
the strong-graded-equivalence $(A,B,X,Y,$ $\tau _A,\tau _B)$ between $A$
and $B$, and let $e\in \mult _1(\mathbb {L})$ be the multiplier defined
in the first part of the proof of
Proposition~\ref {prop:moridalgebras}. Then the pair $\mathbb {L}$ and $e$
satisfy (2) of Proposition~\ref {prop:moridalgebras}. Moreover, using
Proposition~\ref {prop:strong}:
\[D_t^{\mathbb {L}}
=\begin {pmatrix}A_t&X_t\\
  Y_t&B_t\end {pmatrix}\begin {pmatrix}A_{t^{-1}}&X_{t^{-1}}\\
  Y_{t^{-1}}&B_{t^{-1}}\end {pmatrix}
=\begin {pmatrix}D_t^A&D_t^AX_1\\
  D_t^BY_1&D_t^B\end {pmatrix}.
\]
Thus $eD_t^{\mathbb {L}}e\cong D_t^A$, $(1-e)D_t^{\mathbb {L}}(1-e)\cong
D_t^B$.  Now, in order to compute
$D_t^{\mathbb {L}}eD_t^{\mathbb {L}}$,
note first that  it  follows from (2) of Proposition~\ref {prop:strong} that
$D_t^BY_1D_t^A=D_t^BY_1$ and
$\tau _B(D_t^BY_1\otimes _AD_t^AX_1) =D_t^B B_1 = D_t^B$.  Therefore:
\begin {gather*}
D_t^{\mathbb {L}}eD_t^{\mathbb {L}}
=\begin {pmatrix}D_t^A&D_t^AX_1\\
  D_t^BY_1D_t^A&\tau _B(D_t^BY_1\otimes _AD_t^AX_1)\end {pmatrix}
=D_t^{\mathbb {L}}.
\end {gather*}
The equality between $D_t^{\mathbb {L}}(1-e)D_t^{\mathbb {L}}$ and
$D_t^{\mathbb {L}}$ is proved in a similar way. Finally, let us see that
$\mathbb {L}$ is partially-strongly-graded. Since both $A$
and $B$ are partially-strongly-graded,  we have:
\[
D_t^{\mathbb {L}}\mathbb {L}_t
=\begin {pmatrix}A_t+D_t^A\tau _A(X_1\otimes Y_t)&D_t^AX_t+D_t^AX_1B_t\\
  D_t^BY_1A_t+D_t^BY_t&D_t^B\tau _B(Y_1\otimes X_t)+B_t\end {pmatrix}
\]
Now, using Proposition~\ref {prop:strong} and the fact that $X_t$ is  unital as a left $D_t^A$-module,
we have:
\begin {gather*}
A_t+D_t^A\tau _A(X_1\otimes Y_t)=A_t + D_t^AA_1A_t=A_t+A_t=A_t,\\
D_t^AX_t+D_t^AX_1B_t=X_t+D_t^AX_t=X_t,
\end {gather*}
and similarly we see that $D_t^BY_1A_t+D_t^BY_t=Y_t$ and
$D_t^B\tau _B(Y_1\otimes X_t)+B_t=B_t$. Thus
$\mathbb {L}_t\mathbb {L}_t\mathbb {L}_{t^{-1}}=\mathbb {L}_t$.
\par Suppose conversely that $C$ and $e$ satisfy statement (2), and
let $A'=\oplus _{t\in G}A'_t$, $B'=\oplus _{t\in G}B'_t$,
$X=\oplus _{t\in G}X_t$, $Y=\oplus _{t\in G}Y_t$ and $\tau _{A'}$,
$\tau _{B'}$ as in the proof of (2)$\to $(1) in
Proposition~\ref {prop:moridalgebras}. Since $A\cong A'$ and $B\cong
B'$ as graded algebras, then both $A'$ and $B'$ are partially-strongly-graded,
so in particular $D_t^{A'}$ and $D_t^{B'}$ are idempotent. Notice  that
\begin {gather*}
 D_t^{A'} = A'_t A'_{t^{-1}} =  eC_t e C_{t^{-1}}e = eC_t D^C_{t^{-1}} e D^C_{t^{-1}} C_{t^{-1}}e =\\
eC_t D^C_{t^{-1}} C_{t^{-1}}e = e D^C_{t} e.
\end {gather*}
 Let us see
that the graded Morita context $(A',B',X,Y,\tau _{A'},\tau _{B'})$ is
strong. In first place, since $C$ is partially-strongly-graded  and satisfies (2)(b), we have
\begin {gather*}
D_t^{A'}X_t=(eC_t C_{t^{-1}}e)(eC_t(1-e))\\
= eC_tC_{t^{-1}}D_t^CeD_t^CC_t(1-e)
= eC_t C_{t^{-1}}D_t^CC_t(1-e)\\
%=eC_tD_{t^{-1}}^Ce(D_{t^{-1}}^C)^2(1-e)
=eC_t(1-e)
=X_t,
\end {gather*} showing that    $D_t^{A'}X_t=X_t$.
Similar computations yield $X_tD_{t^{-1}}^{B'}=X_t$,
$D_t^{B'}Y_t=Y_t$, and $Y_tD_{t^{-1}}^{A'}=Y_t$.
\par On the other hand:
\begin {gather*}
X_tY_{t^{-1}}=(eC_t(1-e))((1-e)C_{t^{-1}}e)
=eC_tD_{t^{-1}}^C(1-e)D_{t^{-1}}^CC_{t^{-1}}e\\
=eC_tD_{t^{-1}}^CC_{t^{-1}}e
=e(D_t^C)^2e
=eD_t^Ce
=D_t^{A'}.
\end {gather*}
The equality $Y_{t^{-1}}X_t=D_{t^{-1}}^B$ is proved in the same
way. Thus we conclude that $(A',B',X,Y,$ $ \tau _{A'},\tau _{B'})$ is
a strong-graded-equivalence between $A'$ and $B'$, which are
isomorphic to $A$ and $B$ respectively. Hence $A$ and $B$
are strongly-graded-equivalent.
\end {proof}

\begin {prop}\label {prop:sge is eq}
Strong-graded-equivalence is an equivalence relation in the class of partially-strongly-graded algebras.
\end {prop}
\begin {proof}
We have already seen that every partially-strongly-graded algebra is a strongly-graded autoequivalence.
Since the symmetric property of strongly-graded-equivalence is clear,
we prove that
it is transitive. To this end consider strong-graded-equivalences $\mathsf {M}=(A,A',X,Y,\tau _A,\tau _{A'})$ and
$\mathsf {N}=(A',B,X',Y',\tau '_{A'},\tau _B)$ between $A$ and $A'$ and between $A'$ and $B$ respectively.
Let $\mathsf {MN}=(A,B,\bar {X},\bar {Y},\rho _A,\rho _{B})$ be the graded-equivalence constructed in the proof of Proposition~\ref {prop:ge is eq}. We will show
that  $\mathsf {MN}$ is in fact a strong equivalence. First note that  each $\bar {X}_t$ is a unital
$(D_t^A, D_{t^{-1}}^B)$-bimodule
and  each $\bar {Y}_t$ is a unital
$(D_t^B, D_{t^{-1}}^A)$-bimodule. Indeed, using Proposition~\ref {prop:strong} and \eqref {eq:prp3},
 we see that
\begin {gather*}
X_s\otimes _{A'} X'_{s\m t} = X_s\otimes _{A'} D_{s\m t}^{A'}X'_{s\m t}
  = X_s D_{s\m t}^{A'}  \otimes _{A'} X'_{s\m t}\\
 =  X_1 A'_s D_{s\m t}^{A'}  \otimes _{A'} X'_{s\m t}
=  X_1  D_{t}^{A'}  A'_s \otimes _{A'} X'_{s\m t} =     D_{t}^{A}  X_1 A'_s \otimes _{A'} X'_{s\m t}\\
 =     D_{t}^{A}  X_s \otimes _{A'} X'_{s\m t} \subseteq D_{t}^{A}  \bar {X}_t,
\end {gather*} which shows that $\bar {X}_t =  D_{t}^{A}  \bar {X}_t.$ Similarly, $\bar {X}= \bar {X} D_{t\m }^{B},$
$\bar {Y}_t =  D_{t}^{B}  \bar {X}_t,$ and $ \bar {Y}= \bar {Y} D_{t\m }^{A}.$   \\
We have, for $t\in G$:
\begin {gather*}
\rho _A^t(\bar {X}_t\otimes _{D_{t^{-1}}^B}\bar {Y}_{t^{-1}})
=\sum _{r,s\in G}\rho _A\big ((X_r\otimes _{\color {black} A'}X'_{r^{-1}t})\otimes _B (X_{t^{-1}s^{-1}}\otimes _{\color {black} A'}Y_{s})\big )\\
=\sum _{r,s\in G}\tau _A\big ((X_r\otimes _{A'}\tau '_{A'}(X'_{r^{-1}t}\otimes _B Y'_{t^{-1}s^{-1}})Y_s\big ).
\end {gather*}
Now using again  Proposition~\ref {prop:strong} we obtain:
\begin {gather*}
\tau _A\big ((X_r\otimes _{A'}\tau '_{A'}(X'_{r^{-1}t}\otimes _B Y'_{t^{-1}s^{-1}})Y_s\big )=
\tau _A\big (X_1A'_r\otimes _{A'}A'_{r^{-1}t}A'_{t^{-1}s^{-1}}A'_sY_1\big )\\
=\tau _A(X_1\otimes _{A'} Y_1)A_rA_{r^{-1}t}A_{t^{-1}s^{-1}}A_s
=A_1D_r^AD_t^AD_{s^{-1}}^A
=D_r^AD_t^AD_{s^{-1}}^A.
\end {gather*}
Therefore $\rho _A^t(\bar {X}_t\otimes _{D_{t^{-1}}^B}\bar {Y}_{t^{-1}})
=\sum _{r,s\in G}D_r^AD_t^AD_{s^{-1}}^A=D_t^A$, $\forall t\in G$.
Analogous computations show that $\rho _B^t(\bar {Y}_t\otimes _{D_{t^{-1}}^A}\bar {X}_{t^{-1}})=D_t^B$, $\forall t\in G$,
which concludes the proof.
\end {proof}
\par We end this section by showing that the property of having a strong-grading
  is invariant by strongly-graded-equivalence, unlike the situation
  with graded-equivalence (see for
  instance Theorem~\ref {thm:geq} or Theorem~\ref {thm:crossed product
    ppa and env}).

\begin {prop}\label {prop:invsgeq}
Let $A=\oplus _{t\in G}A_t$ and $B=\oplus _{t\in G}B_t$ be partially-strongly-graded algebras that are strongly-graded-equivalent algebras. Then $A$ is strongly-graded if and only if
so is $B$.
\end {prop}
\begin {proof}
Let $(A,B,X,Y,\tau _A,\tau _B)$ be a strong-graded-equivalence. By
Proposition~\ref {prop:strongmoridalgebras}, we may suppose that there
exist a partially-strongly-graded algebra $C=\oplus _{t\in G}C_t$ and
$e=e^2\in \mathfrak {M}_1(C)$, such that $A=eCe$ and $B=(1-e)C(1-e)$
are graded subalgebras of $C$, $X=eC(1-e)$, $Y=(1-e)Ce$,
$eD_t^Ce=D_t^A$, $(1-e)D_t^C(1-e)=D_t^B$,
$D_t^CeD_t^C=D_t^C=D_t^C(1-e)D_t^C$,
and $\tau _A$ and $\tau _B$ are given by the product of $C$. Suppose that
$B$ is strongly-graded. Then $Y_tX_{t^{-1}}=D_t^B=B_1$, $\forall t\in
G$. Now, (1) and  (2) of Proposition~\ref {prop:prp} imply that
$B_tY_{t^{-1}}=Y_1$ and $X_t=X_1B_t$ (recall that each $Y_r$, $X_r$,
are unital $D_r^B$ and $D^B_{r^{-1}}$ modules respectively). Hence
$D_t^A=X_tY_{t^{-1}}=X_1B_tY_{t^{-1}}=X_1Y_1=A_1$. Then
$A_t=D_t^AA_t=A_1A_t=A_t$, $\forall t\in G$, so $A$ is strongly-graded, which ends the proof.
\end {proof}

\section {Product partial actions}\label {sec:ppa}

Let us begin by recalling the usual definition of a partial group action  on
an algebra.
\begin {defn}\label {defn:pa} A   partial action $\alpha $ of a
  group  $G$ on
 an algebra $A$ consists of a family of two-sided ideals $D_t$ in $A$
 $(t \in G)$ and algebra isomorphisms $\alpha _t : D_{t^{-1}} \to
 D_t,$ such that for all $s,t \in G$ the following properties are
 verified:
\begin {enumerate}
  \item $\alpha _1$ is the identity isomorphism $A \to A,$
  \item ${\alpha }_s ({D}_{s^{-1}} \cap {D}_t) = {D}_s \cap {D}_{st},$
  \item ${\alpha }_s({\alpha }_t(x)) = {\alpha }_{st}(x)$, for any $x
    \in {D}_{t^{-1}} \cap {D}_{(st)^{-1}}.$
\end {enumerate}
\end {defn}

We say that a partial action $\alpha $ is {\it idempotent} if each
domain $D_t$ is an idempotent ideal.  In this work we will prefer to
replace in the above definition intersections by products as follows.

\begin {defn}\label {defn:ppa} A  product partial action $\alpha $ of a
  group  $G$ on
 an algebra $A$ consists of a family of two-sided ideals $D_t$ in $A$
 $(t \in G)$ and algebra isomorphisms $\alpha _t : D_{t^{-1}} \to
 D_t,$ such that for all $s,t \in G$ the following properties are
 verified:
\begin {enumerate}
  \item $\alpha _1$ is the identity isomorphism $A \to A,$
  \item $ D_t^2 = D_t$, and ${D}_s   {D}_{t} = {D}_t   {D}_{s},$
  \item ${\alpha }_s ({D}_{s^{-1}}   {D}_t) = {D}_s   {D}_{st},$
  \item ${\alpha }_s({\alpha }_t(x)) = {\alpha }_{st}(x)$, for any $x
    \in {D}_{t^{-1}}   {D}_{(st)^{-1}}.$
\end {enumerate}
\end {defn}

For instance, a partial action whose domains are all unital is a product
  partial action. Other examples of product partial actions are
  the partial actions on C*-algebras, since in this case the product
  of two closed ideals   is equal to their intersection.

\par Notice that an idempotent ideal $I$ in an algebra $A$ is a unital
$A$-bimodule. Hence,
given a product partial action  $\alpha =\{\alpha _t
: D_{t\m } \to
D_t, t\in G\}$ of $G$ on $A,$ each ideal $D_t,$ $t\in G,$ is a unital $A$-bimodule.

The following result shows that, under certain
  circumstances, pro\-duct partial actions can be obtained from global
  actions by restriction. We will see in Theorem~\ref {thm:moritaglob}
  that in fact any product partial action can be obtained, up to
  Morita equivalence, in such a way.
\begin {prop}\label {prop:pparest}
Let $A$ be an idempotent ideal in the algebra $B$, and suppose that
$\beta :G\times B\to B$ is an action such that
$A\beta _t(A)=\beta _t(A)A$, $\forall t\in G$. Define
$D_t:=A\beta _t(A)$. Then $\beta _t(D_{t^{-1}})=D_t$, $\forall t\in G$,
and $\alpha :=(\{\alpha _t\},\{D_t\}_{t\in G})$ is a product partial
action, where $\alpha _t(x)=\beta _t(x)$, $\forall t\in G$ and $x\in
D_{t^{-1}}$.
\end {prop}
\begin {proof}
It is clear that $\beta _t(D_{t\m })=D_t$, and also that each $\beta _t(A)$ is idempotent because so is $A$. Moreover
$A\beta _t(A)=\beta _t(A)A$, so we have $D_t^2=A^2\beta _t(A)^2=D_t$, and
\begin {gather*}
D_sD_t
=A\beta _s(A)A\beta _t(A)
=A\beta _t(\beta _{t\m s}(A)A)A
=A\beta _t(A\beta _{t\m s}(A))A \\
=A\beta _t(A)\beta _s(A)A
=D_tD_s.
\end {gather*}
 On the other hand:
\[\alpha _s(D_{s\m }D_t)
=\beta _s(A\beta _{s\m }(A)A\beta _t(A))
=\beta _s(A)A\beta _s(A)\beta _{st}(A)
=D_sD_{st}.
\]
Finally, condition (4) of Definition~\ref {defn:ppa} is obviously
satisfied.
\end {proof}

The product partial action  $\alpha $ obtained in Proposition~\ref {prop:pparest} will be called the {\it restriction} of the
global action $\beta .$ The action $\beta $ is called a {\it globalization}  of $\alpha $ (as well as any global action
isomorphic to $ \beta $). A globalization $\beta $ of a product partial action $\alpha $ is called {\it minimal} if
$B = \sum _{t\in G} \beta _t (A)$, in which case we have $B^2=B$.   Note that each $\beta _t (A)$ is a unital $B$-bimodule, due to the fact that $A$
is an idempotent ideal in $B.$

In \cite {ades}, a partial action $\alpha :=(\{\alpha _t\},\{D_t\}_{t\in G})$ was called  regular if
$$D_{t_1} \cap D_{t_2} \cap \ldots \cap D_{t_n} = D_{t_1}  D_{t_2} \ldots D_{t_n}, \;\;\;\; \forall t_1, t_2, \ldots , t_n \in G.$$ Evidently
any regular partial action is a product partial action. %\todo {It does not seem to be easy to produce  an example of a  non-regular  product partial action.}
  As it was mentioned in \cite[p. 4961]{ades}, any $C^*$-algebraic partial action is regular as well as any partial action on a von Neumann regular ring. Notice that in a von Neumann regular ring the ideals are idempotent and the intersection of any two of them coincides with their product. The same is true for the ideals of the Jacobian algebra ${\mathbb A_n}$  over a field of cha\-rac\-te\-ris\-tic zero introduced by V. Bavula with respect to the Jacobian Conjecture (see \cite[Theorem 3.1, Corollary 3.10]{Bavula2009}). It follows that any partial action on  ${\mathbb A}_n$ is regular. Moreover, the same can be said about the algebra ${\mathbb I}_n$ of integro-differential operators on a
polynomial algebra thanks to the ideal equivalence with ${\mathbb A_n}$ established in \cite[Theorem 3.1]{Bavula2011}. Furthermore,
\cite[Corollary 7.4]{Bavula2010} implies that any idempotent
partial action on the  algebra of one-sided inverses of a polynomial algebra is regular.

%The next
 % example shows that not every product partial action is necessarily
  %regular:

%\begin {exmp}
%\todo {to do} Not every product partial action is regular.
%\end {exmp}

Given a partial action or a product partial action $\alpha $ of $G$ on $A$, we define the skew group algebra
$A \rtimes _{\alpha }G$
of $A$ by $\alpha $ as the direct sum $\oplus _{t\in G} D_t \delta _t,$ with the product determined  by the rule
\begin {equation}\label {eqn:product}
(a \delta _r)(b \delta _s) = \alpha _r ( {\alpha }\m _r  (a) b)
\delta _{rs},\end {equation}
 where $r,s \in G, a\in D_r, b\in D_s.$
It was established in  \cite [Corollary 3.2]{dokex} that $A \rtimes _{\alpha }G$   is associative if $\alpha $ is an
idempotent partial
action of a group $G$ on a unital algebra $A.$  Notice that the proof does not use the fact that $A$ is unital and,
moreover, it works for any product partial  action  $\alpha .$

Let   $\alpha =\{\alpha _t : D_{t^{-1}} \to
D_t, t\in G\}$  and
${\alpha }' =\{{\alpha }' _t : D'_{t^{-1}} \to D'_t, t\in G\}$ be product partial actions of $G$ on
$A$ and $A',$ respectively. By a {\it morphism} $\varphi : \alpha \to {\alpha }'$
 we mean an algebra homomorphism $\varphi : A \to A'$  such that
 $\varphi (D_t) \subseteq D'_t$  and
$\varphi (\alpha _t (a) ) = {\alpha }'_t (\varphi (a))$ for all $t\in
G,$  $a \in D_{t\m }.$  A morphism  $\varphi : \alpha \to
  {\alpha }'$ induces an algebra homomorphism
  $\varphi ^\rtimes :A\rtimes _\alpha
  G\to A\rtimes _{\alpha '}G$. In fact the correspondence
$(\alpha \stackrel {\varphi }{\to }\alpha ')\mapsto (A\rtimes _\alpha
  G\stackrel {\varphi ^\rtimes }{\to } A\rtimes _{\alpha '}G)$ is a
functor.

\begin {rem}\label {rem:eneveloping ac}
  Every partial action  $\alpha :=(\{\alpha _t\},\{D_t\}_{t\in
    G}),$ such that each   $D_t$ is unital, can be seen as  a restriction of a global
  action, which is essentially unique if a minimality condition is required (see \cite {dokex}).  We
  refer to this global action as the enveloping action of
  $\alpha $. Note that the enveloping algebra is idempotent \cite [Theorem 3.1]{DdRS}, which
  implies that so is the corresponding skew group algebra of the
  enveloping action.
\end {rem}
%\todo {Is a minimal globalization of a ppa necessarily unique?}\todo {I think that this question would be more reasonable
% to consider if we would have a result on the existence of a globalization.  I think that we are ok in this paper, as our main
%spirit here is  Morita equivalence and we  have globalization and uniqueness up to Morita equivalence.}

\begin {prop}\label {prop:ppa-prp}
  The skew group algebra of a product partial action is partially-strongly-graded,  and it is strongly-graded if and only if the
    partial action is global.
\end {prop}
\begin {proof}
 Let $\alpha :=(\{\alpha _t\},\{D_t\}_{t\in G})$ be a product partial
 action on the algebra $A$. Let $B:=A\rtimes _{\alpha }G$. Then
 $B=\oplus _{t\in G}B_t$, where  $B_t=D_t\delta _t$, and the product is
 given by
 $d_s\delta _s\,d_t\delta _t=\alpha _{s}(\alpha _{s}^{-1}(d_s)d_t)\delta _{st}$.
 Thus
 $B_sB_t=\alpha _{s}(\alpha _{s}^{-1}(D_s)D_t)\delta _{st}=\alpha _{s}(D_{s^{-1}}D_t)\delta _{st}=D_sD_{st}\delta _{st}$. In
 particular $B_1B_t=B_t=B_tB_1$, $\forall t\in G$. Moreover we
 have \[B_tB_{t^{-1}}B_t=(D_tD_{\textcolor {black}{1}}\delta _{\textcolor {black}{1}})(D_t\delta _t)=D_t\delta _t=B_t.\]
 As for the last statement, it is clear that the skew group
  algebra of a global action is strongly-graded. To prove the converse
  note that, according to (3) of Definition~\ref {defn:ppa} and
  equality \eqref {eqn:product} above, we have
  $$(D_r\delta _r) (D_s\delta _s)=D_rD_{rs}\delta _{rs}, \quad \forall r,s\in
  G,$$
which implies $D_{rs}\subseteq D_r$, $\forall r,s\in G$, hence
$D_r=A$, $\forall r$.
\end {proof}

Moreover, we have the next.

\begin {prop}\label {prop:pa-prp}
  Let $\alpha $ be an idempotent partial action of a group  $G$ on an algebra $A.$ Then $A \rtimes _{\alpha }G$
is partially-strongly-graded  if and only if
	$\alpha $ is a product partial action.
\end {prop}
\begin {proof} The `if' part follows from Proposition~\ref {prop:ppa-prp}. For the `only if' part suppose that
$B= A\rtimes _{\alpha }G$ satisfies the \prp .  If we write
 $B=\oplus _{t\in G}B_t$, where  $B_t=D_t\delta _t,$ then we are assuming that
  \[B_tB_{t^{-1}}B_t =B_t,\] for all $t\in G.$    Since this is equivalent to the two equalitites in \eqref {eq:prp2}, we obtain,
on  one hand,
	$$B_{r\m } B_r B_s  = D_{r\m } D_s \delta _s = B_{r\m } B_{rs} = \alpha _{r\m } (D_r D_{rs}) \delta _s,$$ which gives
	\begin {equation}\label {eq:pa-prp}
	\alpha _{r} (D_{r\m } D_{s}) =  D_{r} D_{rs}
	\end {equation}for all $r,s \in G.$ On the other hand,
	$$B_r B_s B_{s\m } = \alpha _r (D_{r\m } D_s) \delta _r = B_{rs} B_{s\m } = \alpha _{rs} ( D_{s\m r\m } D_{s\m } ) \delta _r,$$
	that is
	$$\alpha _r (D_{r\m } D_s) =  \alpha _{rs} ( D_{s\m r\m } D_{s\m } ), $$ for all $r,s \in G.$ Thanks to \eqref {eq:pa-prp} the latter results in
	$$ D_r D_{rs} = D_{rs} D_r,$$ for all $r,s \in G,$ showing
        that the domains $D_t$ commute with each other. In combination
        with
	\eqref {eq:pa-prp}  this shows that $\alpha $ is a product partial action. \end {proof}

\begin {thm}\label {thm:crossed product ppa and env}
Let $\alpha $ be a product partial action on $A$ with a minimal
globalization $\beta $ acting on $B$. Then $A\rtimes _\alpha G$ and
$B\rtimes _\beta G$ are graded-equivalent.
\end {thm}
\begin {proof}
Taking the subsets  $X:=\oplus _{t\in G}\beta _t(A)\delta _t  $ and
$Y:=\oplus _{t\in
  G}A\delta _t  $ in $  B\rtimes _\beta G,$ we show that $\big (B\rtimes _\beta G, A\rtimes _\alpha G,
X, Y \big )$  is a graded-equivalence. We only need to show that
it is a Morita equivalence, for all the objects of the context are
graded under the grading of $B\rtimes _\beta G$. Using the fact that each  $\beta _t(A)$ is a unital $B$-bimodule,  we have:
\begin {gather*}
(B\rtimes _\beta G)X
=\sum _{s,t}(B\delta _s)(\beta _t(A)\delta _t)
=\sum _{s,t} B\beta _{st}(A)\delta _{st}
=\sum _{s,t} \beta _{st}(A)\delta _{st}
=X.\\
Y(B\rtimes _\beta G)
=\sum _{s,t}(A\delta _s)(B\delta _t)
=\sum _{s,t}A\beta _s(B)\delta _{st}
=\sum _{t}A\delta _{t}=Y.\\
XY
%=\sum _{s,t}(\beta _s(A)\delta _s)(A\delta _t)
=\sum _{s,t}\beta _s(A)^2\delta _{st}
=\sum _{r}\sum _{t}\beta _{rt\m }(A)\delta _{r}
=\sum _{r}B\delta _{r}
=B\rtimes _\beta G.\\
YX
=\sum _{s,t}(A\beta _{st}(A)\delta _{st}
=\sum _{s,t}D_{st}\delta _{st}
=A\rtimes _\alpha G.
\end {gather*}
Finally:
\begin {gather*}
X(A\rtimes _\alpha G)
=X(YX)=(XY)X=(B\rtimes _\beta G)X=X.\\
(A\rtimes _\alpha G)Y=(YX)Y=Y(XY)=Y(B\rtimes _\beta G)=Y.
\end {gather*}
\end {proof}

  \section {Partial smash products}\label {sec:smash}
\par
Let $B=\oplus _{t\in G}B_t$ be a graded algebra.
It is readily checked that the set $\ps {B}:=\oplus _{r,s\in
  G}B_{r^{-1}}B_s\munit _{r,s}$ is a two-sided ideal of $B\#
G$.

\begin {defn}\label {defn:psmash}
We will say that $\ps {B}:=\oplus _{r,s\in G}B_{r^{-1}}B_s\munit _{r,s}$
is the partial smash product of the $G$-graded algebra $B$. We
may occasionally denote it also by $B\#_pG$.
\end {defn}

Observe that $\ps {B}=B\#G$ if and only if $B$ is strongly-graded.
The easy proof of the next fact is left to the reader:

  \begin {prop}\label {prop:idempotent} Let $B=\oplus _{t\in G}B_t$ be a
    graded algebra, such that   the  $B_1$-bimodule $B_t$ is  unital for each $t\in G.$
Then the algebras $B,$  $B_1$ and   $\ps {B}$ are  idempotent.
% if and only if
  %  $\sum _{s\in G}B_{s^{-1}}B_s=B_1$. In particular $\ps {B}$ is idempotent
 %      whenever so is $B$.
  \end {prop}

  We proceed with the following fact.

  \begin {prop}\label {prop:restsmash}
    Let $B=\oplus _{t\in G}B_t$ be a partially-strongly-graded  algebra and
    $I:=\ps {B}$, its partial smash product. Then
    \begin {enumerate}
      \item The linear $\beta ^B$-orbit of $I$ is all of $B\#G$.
      \item For every
        $t\in G$ we have
\[I
    \beta ^B_t(I)=\beta ^B_t(I)
    I=\oplus _{r,s}B_{r^{-1}}B_tB_{t^{-1}}B_s\munit _{r,s},\]   so that
    the restriction $\gamma ^B:=\beta ^B |_I$ of $\beta ^B $
to $I$ is a product
    partial action, and $\beta ^B$ is a minimal globalization of
    $\gamma ^B$.
    \end {enumerate}
  \end {prop}
  \begin {proof} Since $B$ is partially-strongly-graded, each $B_t$ is a unital $B_1$-module
and  by
    Proposition~\ref {prop:idempotent} the ideal $I$ is
    idempotent. Therefore
\[B\# G=\oplus _{r,s\in G}\beta _s^B(B_{r^{-1}s}B_1e_{s^{-1}r,1})
\subseteq \sum _{s\in G}\beta _s^B(I)
\subseteq B\# G,
\]which
proves our first statement.
\par Let
us see next that $\gamma $ is a product partial action. On  one hand:
\[\beta ^B_t(I)=\oplus _{r,s}B_{r^{-1}}B_s\munit _{tr,ts}=\oplus _{u,v}B_{u^{-1}t}B_{t^{-1}v}\munit _{u,v}.\]
Then, since $ B_1=\sum _{u\in G}B_{u^{-1}}B_u$, we have:
\begin {gather*}
I  \beta _t^B(I)
=\sum _{r,s,u,v}B_{r^{-1}}B_sB_{u^{-1}t}B_{t^{-1}v}\munit _{r,s}\munit _{u,v}
=\sum _{r,u,v}B_{r^{-1}}B_uB_{u^{-1}t}B_{t^{-1}v}\munit _{r,v}\\
=\sum _{r,u,v}B_{r^{-1}}B_uB_{u^{-1}}B_tB_{t^{-1}}B_v\munit _{r,v}
=\oplus _{r,s}B_{r^{-1}}B_tB_{t^{-1}}B_s\munit _{r,s}
\end {gather*}
On the other hand:
\begin {gather*}
 \beta _t^B (I)  I
=\sum _{r,s,u,v}B_{r^{-1}t}B_{t^{-1}s}B_{u^{-1}}B_v\munit _{r,s}\munit _{u,v}
=\sum _{r,u,v}B_{r^{-1}t}B_{t^{-1}u}B_{u^{-1}}B_v\munit _{r,v}\\
=\sum _{r,u,v}B_{r^{-1}}B_tB_{t^{-1}}B_uB_{u^{-1}}B_v\munit _{r,v}
=\oplus _{r,s}B_{r^{-1}}B_tB_{t^{-1}}B_s\munit _{r,s}.
\end {gather*}
Thus $I  \beta ^B_t(I)=  \beta _t^B(I)   I$, so $\gamma ^B$ is a product
partial action by Proposition~\ref {prop:pparest}, and $\beta ^B$ is a
minimal globalization of $\gamma ^B$ by~(1).
  \end {proof}

Note that in the proof of (1) of Proposition~\ref {prop:restsmash} the only restriction on the grading of $B$
we used is that each $B_t$ is a unital $B_1$-bimodule.  Thus we
  have:
\begin {cor}\label {cor:sgchar}
Let $B=\oplus _{t\in G}B_t$ be a graded algebra such that every $B_t$ is
a unital $B_1$-bimodule. Then $B\#_pG$ is $\beta ^B$-invariant if
and only if $B\#_pG=B\#G$, that is, if and only if $B$ is
strongly-graded.
\end {cor}

  \begin {thm}\label {thm:partialrep}
If $B=\oplus _{t\in G}B_t$ is a partially-strongly-graded algebra, then
$(B\#G) \rtimes _{{\beta }^B} G$ and $\ps {B}\rtimes _{{\gamma }^B} G$ are
graded-equivalent.
\end {thm}
  \begin {proof} This is a direct consequence of
    Proposition~\ref {prop:restsmash} and
    Theorem~\ref {thm:crossed product ppa and env}.
  \end {proof}

  We shall call  $\gamma ^B$ from Proposition~\ref {prop:restsmash}
the {\it canonical partial action} of $G$ on the partial smash product
$\ps {B}.$

{If $\phi :A\to B$ is a homomorphism of $G$-graded algebras and
  $\phi ^\#:A\# G\to B\# G$ is the corresponding homomorphism between
  the smash products, it is clear that $\phi ^\#(\ps {A})\subseteq
  \ps {B}$. Thus $\phi ^\#$ induces a homomorphism $\phi ^\natural :\ps {A}\to
  \ps {B}$. Besides, since
  $\beta _t^B\phi ^\#=\phi ^\#\beta _t^A$, we have that
  $\gamma _t^B\phi ^\natural =\phi ^\natural \gamma _t^A$ $\forall t\in G$,
  so that
  $\phi ^\natural :\gamma ^A\to \gamma ^B$ is a homomorphism of partial
  actions. Therefore $\phi ^\natural $ induces a homomorphism
  $\phi ^\rtimes :\ps {A}\rtimes G\to \ps {B}\rtimes G$. It turns out that the
  maps $(A\stackrel {\phi }{\to }B)\mapsto
  (\gamma ^A\stackrel {\phi ^\natural }{\to }\gamma ^B)$ and
  $(A\stackrel {\phi }{\to }B)\mapsto
  (\ps {A}\rtimes G\stackrel {\phi ^\rtimes }{\to }\ps {B}\rtimes G)$ are functors
  from the category of graded algebras to the category of partial
  actions on algebras and to that of
	 partially-strongly-graded
  algebras, respectively.}

\begin {prop}\label {prop:duality}
Let $B=\oplus _{t\in G}B_t$ be a partially-strongly-$G$-graded
algebra. Then $\ps {B}\rtimes _{\gamma ^B}G$ is naturally isomorphic, as
a graded algebra, to the graded subalgebra
$\oplus _{r,s,t}B_{r^{-1}}B_tB_se_{r,s}$ of ${\rm FMat}_G(B)$.
\end {prop}
\begin {proof}
Let $\psi _B:(B\#G)\rtimes _{\beta ^B}G\to {\rm FMat}_G(B)$ be the
natural isomorphism defined in Theorem~\ref {thm:duality}.
Recall that
$\psi _B$ was determined by
$be_{r,s}\delta _t\stackrel {\psi _B}{\to }be_{r,s}\Delta _t=be_{r,t^{-1}s}$, $\forall
r,s,t\in G$, $b\in B_{r^{-1}s}$.
Since $\ps {B}\rtimes _{\gamma ^B}G=\oplus _{t\in G}I_t\delta _t$, where $I_t=I
    \beta ^B_t(I)=\beta ^B_t(I)
    I=\oplus _{r,s}B_{r^{-1}}B_tB_{t^{-1}}B_s\munit _{r,s}$, it
follows that $\psi _B(\ps {B}\rtimes _{\gamma ^B}G)=\oplus _{t\in
  G}\psi _B(I_t\delta _t)=\oplus _{t\in G}I_t\Delta _t$. Now:
\begin {gather*}
I_t\Delta _t
=\oplus _{r,s\in G}B_{r^{-1}}B_tB_{t^{-1}}B_s\munit _{r,t^{-1}s}
=\oplus _{r,u}B_{r^{-1}}B_tB_{t^{-1}}B_{tu}e_{r,u}\\
=\oplus _{r,s}B_{r^{-1}}B_tB_{t^{-1}}B_{t}B_se_{r,s}
=\oplus _{r,s}B_{r^{-1}}B_{t}B_se_{r,s}.
\end {gather*}
To verify naturality, just note that if $\phi :A\to B$ is a
homomorphism of partially-strongly-graded algebras, the natural map
$\tilde {\phi }: (A\#G)\rtimes _{\beta ^A}G\to (B\#G)\rtimes _{\beta ^B}G$
sends $\ps {A}\rtimes _{\gamma ^A}G$ into $\ps {B}\rtimes _{\gamma ^B}G$, so the
following diagram commutes:
%\todo {Correction in the diagram: place's name of the first arrow in
 % the second row}
\[\xymatrix { (A\#_pG)\rtimes _{\gamma ^A}G\ar@ {->}[r]^{inc}\ar [d]_-{\tilde {\phi }}&(A\# G)\rtimes _{\beta ^A}G
\ar@ {->}[r]^-{\psi _A}_-{\cong }\ar [d]_-{\tilde {\phi }}
&{\rm FMat}_G(A)\ar [d]^{\phi ^{\textrm {fin}}}\\
(B\#_pG)\rtimes _{\gamma ^B}G\ar@ {->}[r]_{inc}&(B\# G)\rtimes _{\beta ^B}G\ar@ {->}[r]_-{\psi _B}^-{\cong }
&{\rm FMat}_G(B)} \]
\end {proof}

\section {Equivalences of product partial actions.}\label {MoritaEqParAc}
  \subsection {Morita equivalence.}\label {ssMoritaEqParAc}

If $\theta =\{\bar {D}_{t\m }\stackrel {\theta _t}{\to }\bar {D}_{t}\}_{t \in G} $ is a product
partial action of the group $G$ on the algebra $C$, we say that a
subset $S$ of $C$ is $\theta $-invariant if $\theta _t(S\cap
\bar {D}_{t\m })=S\cap \bar {D}_{t}$, $\forall t\in G$. Suppose that
$\mathsf {M}=(A,A',X,Y,\tau _A,\tau _{A'})$ is a Morita
context. If $X'\subseteq X$ and $Y'\subseteq Y$, in general we will
write $X'Y'$ instead of $\tau _A(X'\otimes _{A'}Y')$ and $Y'X'$ instead
of $\tau _{A'}(Y'\otimes _{A}X')$.

\begin {defn}\label {defn:Moritapactions} Let $\alpha = \{
  {D}_{t\m }\stackrel {{\alpha }_t}{\to }  {D}_t\}_{t \in G}$ and $\alpha ' = \{
  {D'}_{t\m }\stackrel {{\alpha '}_t}{\to }  {D'}_t\}_{t \in G}$ be
  product  partial actions of $G$ on algebras $A$ and ${A
  }',$ respectively.  We say that $\alpha $ is  Morita equivalent to
  ${\alpha }'$  if there exists a Morita equivalence
  $\mathsf {M}=(A,A',X,Y,\tau _A,\tau _{A'})$ between $A$ and $A'$, and a
  product partial
  action $\theta =\{\bar {D}_{t\m }\stackrel {\theta _t}{\to }\bar {D}_t\}$
  on the context algebra $C$ of $\mathsf {M}$ such that:
  \begin {enumerate}
  \item [(i)] $Y{D}_t X = {D}'_t$.
  \item [(ii)] $\theta |_A=\alpha $ and $\theta |_{A'}=\alpha '$.
  \end {enumerate}
The pair $(\mathsf {M},\theta )$ will be called a Morita
equivalence between $\alpha $ and $\alpha '$. We simbolize this relation
by writing $\alpha \stackrel {M}{\sim }\alpha '$
\end {defn}

\begin {rem}\label {rem:note that item}
Note that the above definition implies that
%the algebras $A $ and
%${A}'$ are Morita equivalent.
each  ideal ${D}_t$ corresponds to
${D}'_t$ under this equivalence. In fact, condition (i) above implies
\[D_t=(XY)D_t(XY)=X(YD_tX)Y=X {D}'_t Y.\] Moreover, \[{D}_t X= D_tX{D}'_t=X{D}'_t\] and
\[{D}'_t Y = {D}'_t YD_t=Y{D}_t,\quad {} \forall t\in G.\] Define
$X_t=D_tX{D}'_t$ and $Y_t={D}'_t YD_t$. Then
\begin {gather}
X_tY_t=D_t \qquad \textrm { and }\qquad Y_tX_t=D'_t, \quad \forall t\in G. \label {eq:Moritapactions1}\\
D_tX=X_t=XD'_t\qquad \textrm { and }\qquad D'_tY=Y_t=YD_t,\quad \forall t\in G.\label {eq:Moritapactions2}
\end {gather}

 \begin {rem}
 It follows from the previous remark that if $\alpha $ and $\alpha '$ are Morita equivalent
 product partial actions, then $\alpha $ is a global
 action if and only if so is $\alpha '$.
 \end {rem}

\begin {rem}\label {rem:it is easy}
It is easy to see, using Remark~\ref {rem:note that item}, that in Definition~\ref {defn:Moritapactions}
condition (i) is equivalent to the following one:

\end {rem}
(i') There exist families $\{X_t\}_{t\in G}$ and $\{Y_t\}_{t\in G}$ of
subbimodules of $X$ and $Y$ respectively such that
\eqref {eq:Moritapactions1} and \eqref {eq:Moritapactions2} hold.

\end {rem}
\begin {rem}\label {rem:global}
Observe that the proof of
 \cite [Proposition~2.11]{ades}  works for   product partial actions. It  follows, in particular, that
  $\bar {D}_t=\begin {pmatrix}D_t&X_t\\
     Y_t&D_t'\end {pmatrix}$, $\forall t\in G$, and all of the subsets
  $\begin {pmatrix} A& 0\\ 0& 0\end {pmatrix}$, $\begin {pmatrix} 0&0 \\
   0 & A'\end {pmatrix}$, $\begin {pmatrix} 0& X\\0 &0 \end {pmatrix}$ and
  $\begin {pmatrix}0 & 0\\Y & 0\end {pmatrix}$ of the context algebra
  $C$ are $\theta $-invariant. We will identify the subsets above with
  $A$, $A'$, $X$ and $Y$ respectively. Notice that with this identification and
 slightly abusing the notation we may write
\begin {equation}\label {alpha-module}
\theta _t (ax)= \alpha _t(a) \, \theta _t (x),  \theta _t (x a')=  \theta _t (x) \, {\alpha }' _t(a'),
\end {equation} for all $a\in D_t^A, x \in X, a' \in D_t^{A'},$ and similar equalities can be written for the
bimodule $Y.$ Observe also that if
  both of $\alpha $ and $\alpha '$ are global actions, then $\theta $
  must be a global action as well.
\end {rem}

Let $(\mathsf {M},\theta )$ be a Morita equivalence between the
 product
  partial actions $\alpha $ and $\alpha '$, where $\mathsf {M}=(A, A', X,
  Y, \tau _A, \tau _{A'})$. In what follows we use the notation of
  Remark~\ref {rem:note that item}. Every $X_t$ is a unital
  $(D_t,D'_t)$-bimodule, and every $Y_t$ is a unital
  $(D'_t,D_t)$-bimodule. Therefore we have natural maps
  $j_t:X_t\otimes _{D'_t}Y_t\to X\otimes _{A'}Y$ and
  $j'_t:Y_t\otimes _{D_t}X_t\to Y\otimes _{A}X$. We define
  $\tau _t:=\tau _Aj_t:X_t\otimes _{D'_t}Y_t\to D_t$ and
  $\tau _t':=\tau _{A'}j_t':Y_t\otimes _{D_t}X_t\to D_t'$. Then $\tau _t$
  is a $D_t$-bimodule map, and $\tau _t'$ is a $D_t'$-bimodule
  map. Let us see that $\tau _t$ is surjective. Since $\tau _A$ is
  surjective, given $a\in A$ there exist $x_1,\ldots ,x_n\in X$ and
  $y_1,\ldots ,y_n\in Y$ such that
  $\tau _A(\sum _{i=1}^nx_i\otimes y_i)=a$. Now, given $d_1,d_2\in D_t$,
we have $d_1x_1,\ldots ,d_1x_n\in X_t$, $y_1d_2,\ldots ,y_nd_2\in Y_t$,
and $\tau _t(\sum _{i=1}^nd_1x_i\otimes y_id_2)=d_1ad_2$. Thus
$\tau _t(X_t\otimes _{D_t'}Y_t)\supseteq D_tAD_t=D_t$, because $D_t$ is
idempotent and unital over $A$. Similarly, every $\tau _t'$ is surjective. Then we have:

\begin {prop}\label {prop:moreqideals}
The context
$\mathsf {M}_t:=(D_t,D_t',X_t,Y_t,\tau _t,\tau _t')$ is a
Morita equivalence between $D_t$ and $D_t'$.
\end {prop}

\begin {rem}\label {def2MoritaEq} Observe that the alternative definition of the Morita equi\-va\-len\-ce of regular
partial actions
given in    \cite [Proposition~ 2.11]{ades} also holds  for product partial actions without any change in the
proof.
\end {rem}

\begin {prop}\label {prop:equivalencerel}   Morita equivalence of
  product partial actions  is an equivalence relation.
\end {prop}
\begin {proof}
  Just follow, \textit {mutatis mutandis}, the proof in
   \cite [Proposition~ 2.12]{ades} of the same result for regular
   partial actions, taking into account Remark~\ref {def2MoritaEq}.
  \end {proof}

 \begin {prop}\label {prop:moritaeqs}
 If $\alpha $ and $\alpha '$ are Morita equivalent product partial
 actions on $A$ and $A'$ respectively, then $A\rtimes _\alpha G$ and
 $A'\rtimes _{\alpha '}G$ are strongly-graded-equivalent (the converse
 is also true: see Theorem~\ref {thm:moritaskew}).
 \end {prop}
 \begin {proof}
We use the notation of Definition~\ref {defn:Moritapactions}. Let
 $\bar {A}:=A\rtimes _\alpha G$, $\bar {B}:=A'\rtimes _{\alpha '} G=$,
 and $\bar {C}:=C\rtimes _\theta G$. Thus $\bar {A}=\oplus _{t\in
   G}D_t\delta _t$, $\bar {B}=\oplus _{t\in G}D_t'\delta _t$, and
 $\bar {C}=\oplus _{t\in G}\bar {D}_t\delta _t $.   Consider
 $e:=(L,R)\in \mult _1(\bar {C})$,
 where $L\big (\begin {pmatrix}a&x\\ y&b\end {pmatrix}\delta _t\big )
=\begin {pmatrix}a&x\\0&0\end {pmatrix} \delta _t$ and
$R\big (\begin {pmatrix}a&x\\y &b\end {pmatrix}\delta _t\big )
=\begin {pmatrix}a&0\\y &0\end {pmatrix}\delta _t$.
To see that $\bar {A}$ and $\bar {B}$ are strongly-graded-equivalent, it
is enough to see that the pair $(\bar {C},e)$ satisfies (2) of
Proposition~\ref {prop:strongmoridalgebras}. First note that
\begin {gather*}
\bar {C}_re\bar {C}_s
=\begin {pmatrix}D_r&X_r\\Y _r&D'_r\end {pmatrix}\delta _re
\begin {pmatrix}D_s&X_s\\Y _s&D'_s\end {pmatrix}\delta _s
=\begin {pmatrix}D_r&X_r\\Y _r&D'_r\end {pmatrix}\delta _r
\begin {pmatrix}D_s&X_s\\0&0\end {pmatrix}\delta _s\\
=\theta _r\left (\begin {pmatrix}D_{r^{-1}}&X_{r^{-1}}\\Y _{r^{-1}}&D'_{r^{-1}}\end {pmatrix}\begin {pmatrix}D_s&X_s\\0&0\end {pmatrix}\right )\delta _{rs}
=\theta _r\left (\begin {pmatrix}D_{r^{-1}}D_s&D_{r^{-1}}X_s\\Y _{r^{-1}}D_s&Y_{r^{-1}}X_s\end {pmatrix}\right )\delta _{rs}\\
=\theta _r\left (\begin {pmatrix}D_{r^{-1}}D_s&D_{r^{-1}}D_sX\\YD _{r^{-1}}D_s&YD_{r^{-1}}D_sX\end {pmatrix}\right )\delta _{rs}
=\begin {pmatrix}D_{r}D_{rs}&D_{r}D_{rs}X\\YD _{r}D_{rs}&YD_{r}D_{rs}X\end {pmatrix}\delta _{rs}.
\end {gather*}
By \eqref {eq:Moritapactions1}, \eqref {eq:Moritapactions2} and the
computations above we obtain, for $r=t$, $s=1$:
\[\bar {C}_te\bar {C}_1
=\begin {pmatrix}D_{t}D_{t}&D_{t}D_{t}X\\YD _{t}D_{t}&YD_{t}D_{t}X\end {pmatrix}\delta _{t}
%=\begin {pmatrix}D_{t}&D_{t}X_t\\Y _tD_{t}&Y_tX_t\end {pmatrix}\delta _{t}
=\begin {pmatrix}D_{t}&X_t\\Y _t&D'_t\end {pmatrix}\delta _{t}
=\bar {C}_t.
\]
Then it follows that $\bar {C}e\bar {C}=\bar {C}$. Similar computations
show that also $\bar {C}(1-e)\bar {C}=\bar {C}$. On the other hand:
\begin {gather*}
e\begin {pmatrix}a&x\\y &b\end {pmatrix}\delta _t e
=\begin {pmatrix}1&0\\0&0\end {pmatrix}
              \begin {pmatrix}a&x\\y &b\end {pmatrix}
              \begin {pmatrix}1&0\\0&0\end {pmatrix} \delta _t\\
=\begin {pmatrix}a&x\\0&0\end {pmatrix}
 \begin {pmatrix}1&0\\0&0\end {pmatrix} \delta _t
=\begin {pmatrix}a&0\\0&0\end {pmatrix}\delta _t.
\end {gather*}
Similarly, $(1-e)\begin {pmatrix}a&x\\y &b\end {pmatrix}\delta _t (1-e)=
\begin {pmatrix}0&0\\0&b\end {pmatrix}\delta _t$. Then
$e\bar {C}e\cong \bar {A}$ and $(1-e)\bar {C}(1-e)\cong \bar {B}$. Now,
since $D_t^{\bar {C}}=\bar {D}_t\delta _1$, we have
$eD_t^{\bar {C}}e=\begin {pmatrix}D_t&0\\0&0\end {pmatrix}\delta _1\cong
D_t\delta _1=D_t^{\bar {A}}$, and also
 $(1-e)D_t^{\bar {C}}(1-e)\cong
D_t^{\bar {B}}$. Finally, using again the identities
\eqref {eq:Moritapactions1}, \eqref {eq:Moritapactions2} we get:
\begin {gather*}
D_t^{\bar {C}}eD_t^{\bar {C}}
=\bar {D}_t\delta _1e\bar {D}_t\delta _1
=\begin {pmatrix}D_t&X_t\\Y _t&D'_t\end {pmatrix}
\begin {pmatrix}D_t&X_t\\0&0\end {pmatrix}\delta _1
=\begin {pmatrix}D_t&X_t\\Y _t&D'_t\end {pmatrix}\delta _1
=D_t^{\bar {C}},
\end {gather*}
and similarly $D_t^{\bar {C}}eD_t^{\bar {C}}=D_t^{\bar {C}}$, which ends
the proof.
 \end {proof}

 \begin {thm}\label {thm:sg}
Let $B=\oplus _{t\in G}B_t$ be a partially-strongly-graded algebra and
$\ps {B}$ its partial smash product. Then $B$
is strongly-graded-equivalent to the partial skew group algebra
$\ps {B}\rtimes _{{\gamma ^B}} G$.
 \end {thm}
 \begin {proof}
     Let $X_t=\oplus _{r\in G}X_t(r)$, $Y_t=\oplus _{r\in G}Y_t(r)$ and $B':=\oplus _{t\in G}B_t\munit _{1,t}\delta _t$ be the
subsets of $B\# G$ defined in the proof of  Theorem~\ref {thm:geq}, so $X_t(r)=B_{r^{-1}t}\munit _{r,t}$, and
$Y_t(r)=B_r\munit _{1,r}$ $\forall r,t\in G$. As it is shown in Theorem~\ref {thm:geq}, if $X=\oplus _{t\in G}X_t\delta _t$ and
$Y=\oplus _{t\in G}Y_t\delta _t$, then
		$((B\#G) \rtimes G, B',X,Y)$ is a graded-equivalence. We define
		$$X'= (\ps {B}\rtimes _{\gamma ^B} G) X = \oplus _{t\in G} X'_t \delta _t \quad \text {and} \quad Y' = Y (\ps {B}\rtimes _{{\gamma ^B}} G )  = \oplus _{t\in G} Y'_t \delta _t.$$ We are going to show that
$({\ps {B}}\rtimes _{{\gamma ^B}} G, B', X' , Y' )$ gives a strong-graded-equivalence.

We start by recalling   that
the partial action $\gamma {:=\gamma ^B} =\{\gamma _t : I_{t\m } \to I_t \} $ of $G$ on $I = \ps {B}$ is the restriction of $\beta $ with
$$ I_t = \oplus _{r,s} B_{r\m } B_t B_{t\m } B_s \munit _{r,s}.$$  We compute first $X'$ as follows. Notice that
\begin {gather*} (B_{r\m }D^B_u   B_s \munit _{r,s} \delta _u)(B_{v\m u\m
    t}\munit _{v,u\m t} \delta _{u\m t})
=B_{r\m }D^B_u B_s B_{v\m u\m t}\munit _{r,s}\munit _{uv,t}\delta _t\\
= [s=uv] B_{r\m }D^B_u B_s B_{s\m t}\munit _{r,t} \delta _t.
\end {gather*}  Consequently, since
$ X' _t \delta _t = \sum _{u,v} (I_u \delta _u) ( X_{u\m t}(v) \delta
_{u\m t}),$ we obtain
\begin {gather*} X'_t = \sum _{r,s,v} \sum _{u} [s=uv] B_{r\m } (D_u ^B B_s) B_{s\m t} \munit _{r,t}\\
=\sum _{r}  B_{r\m }   \sum _{s}B_s B_{s\m t} \munit _{r,t}=
\sum _{r}  B_{r\m }   B_{ t} \munit _{r,t},
\end {gather*} using Proposition~\ref {prop:idga}. Then $X'=\oplus _{t,r}B_{r\m }B_t\munit _{r,t}\delta _t$. Computing
next $Y',$ we have
\begin {gather*} (B_r \munit _{1,r} \delta _u) ( B_{v\m } D_{u\m t} ^B B_w \munit _{v,w} \delta _{u\m t}) =
B_r \munit _{1,r} \beta _{u} ( B_{v\m } D_{u\m t}^B B_w \munit _{v,w}) \delta _t\\
= [r=uv] B_r  B_{r\m u} D_{u\m t}^B B_w \munit _{1,uw}    \delta _t.
\end {gather*} Setting $s= uw$ we obtain
\begin {gather*}
Y'_t = \sum _{r,s,u,v}  [r=uv] B_r  B_{r\m u} D_{u\m t}^B B_{u\m s} \munit _{1,s}\\
 = \sum _{u,s}(\sum _r B_r B_{r\m u}) D_{u\m t}^B B_{u\m s} \munit _{1,s} =
\sum _{u,s}  B_{ u} D_{u\m t}^B B_{u\m s} \munit _{1,s}\\
=   D_{ t}^B \sum _{s} (\sum _u  B_{ u} B_{u\m s}) \munit _{1,s} =   D_{ t}^B \sum _{s}   B_{ s} \munit _{1,s},
\end {gather*} using \eqref {eq:prp3}.
%the formula $ B_u D_v^B = D_{uv}^B B_u,$ which is a consequence of the partial representation property.
Thus $Y' = \oplus _{s, t}      D_{ t}^B    B_{ s} \munit _{1,s}  \delta _t.$

Let us see that $X'Y'=I\rtimes _\gamma G$. We have
$$ X'Y'=(I\rtimes _\gamma G)(XY)(I\rtimes _\gamma G)
=(I\rtimes _\gamma G)((B\# G)\rtimes _\beta G)(I\rtimes _\gamma G).$$
Then
%$X'Y'=\oplus _t\sum _{uvw=t}(I_u\delta _u)((B\#G)\delta _v)(I_w\delta _w)$. Hence
\begin {gather*}
X'Y'=\oplus _t\sum _{uvw=t}(I_u\delta _u)((B\#
G)\delta _v)(I_w\delta _w)\\
=\oplus _t\sum _{uvw=t}(I_u\delta _u)((B\# G)\delta _v)(I_w\delta _w)
=\oplus _t\sum _{uvw=t}(I\beta _u(I)(B\# G)\delta _{uv})(\beta _w(I)I\delta _w)\\
=\oplus _t\sum _{uvw=t}I\beta _u(I)\beta _{t}(I)\beta _{uv}(I)\delta _{t}
=\oplus _t\sum _{uvw=t}I_uI_{uv}I_t\delta _t=\oplus _tI_t\delta _t
=I\rtimes _\gamma G.
\end {gather*}

To prove that $Y'X'=B'$, we compute, for each $t,u,r,s\in G$:
\begin {gather*}
(D_u^BB_s\munit _{1,s}\delta _u)(B_{r\m }B_{u\m
  t}\munit _{r, u\m t}\delta _{u\m t})
=[s=ur]D_u^BB_{ur}B_{r\m }B_{u\m t}\munit _{1,t}\delta _t
\end {gather*}
Then, since $D_u^BB_{ur}B_{r\m }=B_uD_r^B$, we have:
\begin {gather*}
Y'X'=\sum _{t,u,r,s}[s=ur]B_uD_r^BB_{u\m t}\munit _{1,t}\delta _t
=\sum _{t,u,r}B_uD_r^BB_{u\m
  t}\munit _{1,t}\delta _t\\
=\sum _{t,u}B_uB_{u\m t}\munit _{1,t}\delta _t
=\oplus _tB_t\munit _{1,t}\delta _t
=B'.
\end {gather*}

Note that $X'$ and $Y'$ are unital over $I\rtimes _\gamma G$ and
$B'$. In fact:
\begin {gather*}
X'B'=X'Y'X'=(I\rtimes _\gamma G)X'
=(I\rtimes _\gamma G)^2X
=(I\rtimes _\gamma G)X =X'.\\
B'Y'=Y'X'Y'=Y'(I\rtimes _\gamma G)
=Y(I\rtimes _\gamma G)^2
=Y(I\rtimes _\gamma G)
=Y'.
\end {gather*}
We have proved so far that $(I\rtimes _\gamma G,B',X',Y')$ is a graded-equivalence. Next we check that it is also a strong-graded-equivalence. On one hand:
\begin {gather*}
(X'_t\delta _t)(Y'_{t^{-1}}\delta _{t\m })
=X'_t\beta _t(Y'_{t^{-1}})\delta _1
=\sum _{r,s}B_{r\m }B_tD_{t\m }^BB_s\munit _{r,t}\munit _{t,ts}\delta _1\\
=\sum _{r,s}B_{r\m }B_tB_s\munit _{r,ts}\delta _1
=\sum _{r,s'}B_{r\m }B_tB_{t\m s'}\munit _{r,s'}\delta _1
=\sum _{r,s}B_{r\m }D_t^BB_s\munit _{r,s}\delta _1\\
=I_t\delta _1 =I_t\beta _{t}(I_{t\m })\delta _1
=(I_t\delta _t)(I_{t\m }\delta _{t\m }).
\end {gather*}
On the other hand:
\begin {gather*}
(Y'_t\delta _t)(X'_{t^{-1}}\delta _{t\m })
=Y'_t\beta _t(X'_{t^{-1}})\delta _1
=\sum _{r,s}D_{t}^BB_rB_{s\m }B_{t\m }\munit _{1,r}\munit _{ts,1}\delta _1\\
=D_{t}^B\sum _{r,s}[r=ts]B_rB_{r\m t}B_{t\m }\munit _{1,1}\delta _1
=D_{t}^B\sum _{r}D_r^B\munit _{1,1}\delta _1\\
=D_{t}^B\munit _{1,1}\delta _1
=(B_t\munit _{1,t}\delta _t)(B_{t\m }\munit _{1,t\m }\delta _{t\m }).
\end {gather*}
\end {proof}

\begin {cor}\label {cor:sg}
If $B=\oplus _{t\in G}B_t$ is a strongly-graded algebra, then $B$ is
strongly-graded-equivalent to $(B\# G){\color {black} \rtimes }_{\beta ^B}G$.
\end {cor}

 By  \cite [(15)]{des} the crossed product by a twisted partial group
 action is partially-strongly-graded.  Thus we have:
\begin {cor}\label {cor:tcp}
 The crossed product by any twisted partial group action is strongly-graded-equivalent to the skew group algebra of a product partial
 action.
    \end {cor}

\begin {thm}\label {thm:moritaE}
Let $\mathsf {M}=(A, B, X, Y, \tau _A, \tau _B)$
be a graded Morita context between $G$-graded algebras $A$ and
$B$. Let $C=\oplus _{t\in G}C_t$ be the corresponding graded context
algebra, so $C_t=\begin {pmatrix} A_t  & X_t \\ Y_t&
  B_t \end {pmatrix}$  Then:
\begin {enumerate}
\item The algebras $C\#G$ and
$\begin {pmatrix}
A\#G  & X\#G \\ Y\#G& B\#G \end {pmatrix}$ are isomorphic,
and the restrictions of $\beta ^C$ to $A\#G$ and $B\#G$ are
  $\beta ^A$ and $\beta ^{B}$ respectively, and $X\#G$, $Y\#G$, are
  $\beta ^C$-invariant.
\item $\mathsf {M}\#G:=(A\#G,B\#G,X\#G,Y\#G)$  is
  a Morita context, where the traces are given by the product in $C\#
  G$. If $\mathsf {M}$ is a Morita equivalence so is $\mathsf {M}\#G$.
\item If $\mathsf {M}$ is a graded-equivalence, then $\beta ^A$ and
  $\beta ^{B}$ are Morita equivalent actions, so
  $(A\#G)\rtimes _{\beta ^A}G$ and $(B\#G)\rtimes _{\beta ^B}G$ are
   strongly-graded-equivalent.
\end {enumerate}
\end {thm}
\begin {proof}
Since $C_{r\m s}=\begin {pmatrix} A_{r\m s}  & X_{r\m s} \\ Y_{r\m s}& B_{r\m
  s} \end {pmatrix}$, it is clear that the linear map $C\#G\to \begin {pmatrix}
A\#G  & X\#G \\ Y\#G& B\#G \end {pmatrix}$ determined by
 \begin {gather*}
c \munit _{r,s}=\begin {pmatrix} a  & x  \\
  y & b  \end {pmatrix}\munit _{r,s}\longmapsto
\begin {pmatrix} a \munit _{r,s}  & x \munit _{r,s} \\
  y \munit _{r,s}& b \munit _{r,s} \end {pmatrix},
\end {gather*} $c \in C_{r\m s},$  is an isomorphism. This proves the first claim of (1). The second
assertion of (1) is clear. Since (3) follows directly
from (1), (2) and Proposition~\ref {prop:moritaeqs},
%\cite [\sout {Theorem 3.1}]{ades}},
%together with the naturality of the canonical
%action on smash products,
we concentrate on the proof of (2). We will suppose that $\mathsf {M}$ is
a Morita equivalence and prove that then so is $\mathsf {M}\# G$.
%It is left to reader to check that our computations below show that
%$\mathsf {M}\# G$ is a Morita context if $\mathsf {M}$ is
%supposed to be just a graded Morita context.
Note that, since $AX=X$ and
$YA=Y$, for each $t\in G$ we have $X_t=\sum _{u}A_uX_{u\m t}$ and
$Y_t=\sum _uY_uA_{u\m t},$ due to the fact that $X$ and $Y$ are graded modules.  Hence
\begin {gather*}
(A\# G)(X\# G)
=\sum _{r,s,u,v}A_{r\m s}\munit _{r,s}X_{u\m v}\munit _{u,v}
=\sum _{r,s,u,v}[s=u]A_{r\m s}X_{s\m v}\munit _{r,v}\\
=\sum _{r,s,v}A_{r\m s}X_{s\m v}\munit _{r,v}
=\sum _{r,u,v}A_{u}X_{u\m r\m v}\munit _{r,v}
=\oplus _{r,v}X_{ r\m v}\munit _{r,v}=X\# G,
\end {gather*}  which shows that  $X\# G$ is a unital left $A\# G$-module.
Similarly  we obtain that $Y\# G$ is a unital right $B\# G$-module:
\begin {gather*}
(Y\# G)(A\# G)
%=\sum _{r,s,u,v}Y_{r\m s}\munit _{r,s}A_{u\m v}\munit _{u,v}
%=\sum _{r,s,u,v}[s=u]Y_{r\m s}A_{s\m v}\munit _{r,v}\\
%=\sum _{r,s,v}Y_{r\m s}A_{s\m v}\munit _{r,v}
=\sum _{r,u,v}Y_{u}A_{u\m r\m v}\munit _{r,v}
=\oplus _{r,v}Y_{ r\m v}\munit _{r,v}=Y\# G.
\end {gather*}
%The proofs of $(B\# G)(Y\# G)=(Y\# G)$ and $(X\# G)(B\# G)=(X\# G)$
%are analogous.
Let us see that the traces are surjective. Recall that $\tau _A$ and
$\tau _B$ are surjective and respect the gradings,  so that
$\tau _A(\sum _{u}(X_{u}\otimes _B Y_{u\m t})=A_t$ and
$\tau _B(\sum _{u}(Y_{u}\otimes _A X_{u\m t})=B_t$, $\forall t\in G$.
Then:
\begin {gather*}
(X\# G)(Y\# G)
=\sum _{r,s,u,v}X_{r\m s}\munit _{r,s}Y_{u\m v}\munit _{u,v}\\
=\sum _{r,s,u,v}[s=u]\tau _A(X_{r\m s}\otimes _BY_{s\m v})\munit _{r,v}
=\oplus _{r,v}A_{r\m v}\munit _{r,v}=A\# G.
\end {gather*}
Similar computations show that $(Y\# G)(X\# G)=B\# G$.

Finally:
\begin {gather*}
(B\# G)(Y\# G)=(Y\# G)(X\# G)(Y\# G)=(Y\# G)(A\# G)=(Y\# G). \\
(X\# G)(B\# G)=( X\# G)(Y\# G)(X\# G)=(A\# G)(X\# G)=(X\# G).
\end {gather*}
\end {proof}

  \begin {prop}\label {prop:reseqppa}
    Let $(\mathsf {M},\gamma )$ be a Morita
    equivalence between the (global) actions $\beta $ and $\beta '$,
    where $\mathsf {M}=(B, B', X, Y, \tau _B, \tau _{B'})$. Let $A$ be an
    idempotent ideal of $B$ such that $A\beta _t(A)=\beta _t(A)A$,
    $\forall t\in G$, and let $\alpha :=\beta |_A$ be the product
    partial action obtained by the
    restriction of $\beta $ to $A$. Let $A':=YAX$ be the ideal that
    corresponds to $A$ via the equivalence $\mathsf {M}$. Then
    \begin {enumerate}
    \item $\beta _s'(A')\beta _t'(A')=Y\beta _s(A)\beta _t(A)X$, $\forall
      s,t\in G$.
    \item Let $\mathsf {N}=(A,A',X_1,Y_1,\tau _A,\tau _{A'})$, where
      $X_1:=AX$, $Y_1=YA$, $\tau _A$ is the composition of $\tau _B$
      with the natural map
      $X_1\otimes _{A'}Y_1\to X\otimes _{B'}Y$, and $\tau _{A'}$ is
      defined analogously. Then $\mathsf {N}$ is a Morita equivalence
      between $A$ and $A'$, and the corresponding context algebra
      $C_{\mathsf {N}}$ is an idempotent ideal of the context algebra
      $C_{\mathsf {M}}$ of $\mathsf {M}$, such that
      $\gamma _s(C_{\mathsf {N}})\gamma _t(C_{\mathsf {N}})=\gamma _t(C_{\mathsf {N}})\gamma _s(C_{\mathsf {N}})$,
      $\forall s,t\in G$.
    \item Let $\alpha ':=\beta '|_{A'}$ and
      $\theta :=\gamma |_{C_{\mathsf {N}}}$ be the restrictions of $\beta '$
      and $\gamma $ to $A'$ and $C_{\mathsf {N}}$ respectively. Then
      $(\mathsf {N},\theta )$ is a Morita equivalence between $\alpha $
      and $\alpha '$.
    \end {enumerate}
\end {prop}
\begin {proof}
Since $X$ and $Y$ are $\gamma $-invariant we have:
\[\beta _t'(A')=\gamma _t(YAX)
=\gamma _t(Y)\gamma _t(A)\gamma _t(X)
=Y\beta _t(A)X.\]
Next observe that every $\beta _t(A)$ is idempotent because so is
$A$. Then
  $B\beta _t(A)=\beta _t(A)$ and, since $XY=B$:
\begin {gather*}Y\beta _s(A)\beta _t(A)X
=Y\beta _s(A)(XY\beta _t(A))X
=(Y\beta _s(A)X)(Y\beta _t(A)X)\\
=\beta _s'(A')\beta _t'(A').
\end {gather*}

 We see next that $\mathsf {N}$ is a Morita equivalence. Since $A$
  and $A'$ are idempotent, then $X_1$ and $Y_1$ respectively a left
  unital $A$-module and a right unital $A'$-module. On the other hand:
\[X_1A'=AX(YAX)=ABAX=AX=X_1, \textrm { and similarly }A'Y_1=Y_1.\]
Reasoning exactly as in the proof of
Proposition~\ref {prop:moreqideals}, we have that $\tau _A$ and
$\tau _{A'}$ are surjective. Let us show that $C_{\mathsf {N}}$ is an
ideal in $C_{\mathsf {M}}$:
\begin {gather*}
\begin {pmatrix}A&AXA'\\
  A'YA&A'\end {pmatrix} \begin {pmatrix}B&X\\Y &B'\end {pmatrix}
=\begin {pmatrix}AB+AXA'Y&AX+AXA'B'\\
  A'YAB+A'Y&A'YAX+A'B'\end {pmatrix}\\
=\begin {pmatrix}A+AXY&AX+AXA'\\
  A'YA+A'Y& A' A'+A'\end {pmatrix}
=\begin {pmatrix}A&AXA'\\
  A'YA&A'\end {pmatrix}
=C_{\mathsf {N}}.
\end {gather*}
A similar computation shows that also
$C_{\mathsf {M}}C_{\mathsf {N}}=C_{\mathsf {N}}$, so $C_{\mathsf {N}}$ is
an ideal of $C_{\mathsf {M}}$. Next we compute
$\gamma _s(C_{\mathsf {N}})\gamma _t(C_{\mathsf {N}})$. Note that (1), with $s=t,$
implies that $\beta _t(A)X=X\beta _t'(A')$ and
$\beta _t'(A')Y=Y\beta _t(A)$, $\forall t\in G$.
Recall also from the proof of
Proposition~\ref {prop:pparest} that, since $A\beta _t(A)=\beta _t(A)A$
$\forall t$, then $\beta _s(A)\beta _t(A)=\beta _t(A)\beta _s(A)$ $\forall
s,t$. Then also $\beta _s'(A')\beta _t'(A')=\beta _t'(A')\beta _s'(A')$ $\forall
s,t$, by (1).
Then  using \eqref {alpha-module}  we have:
\begin {gather*}
\gamma _s(C_{\mathsf {N}})\gamma _t(C_{\mathsf {N}})
=\begin {pmatrix}\beta _s(A)&\beta _s(A)X\\
  \beta _s'(A')Y&\beta _s'(A')\end {pmatrix} \begin {pmatrix}\beta _t(A)&\beta _t(A)X\\
  \beta _t'(A')Y&\beta _t'(A')\end {pmatrix}\\
=\begin {pmatrix}\beta _s(A)\beta _t(A)+\beta _s(A)X\beta _t'(A')Y&\beta _s(A)\beta _t(A)X+\beta _s(A)X\beta _t'(A')\\
  \beta _s'(A')Y\beta _t(A)+\beta _s'(A')\beta _t'(A')Y&\beta _s'(A')Y\beta _t(A)X+\beta _s'(A')\beta _t'(A')\end {pmatrix}\\
=\begin {pmatrix}\beta _s(A)\beta _t(A)&\beta _s(A)\beta _t(A)X\\
  \beta _s'(A')\beta _t'(A')Y&\beta _s'(A')\beta _t'(A')\end {pmatrix}
=\gamma _t(C_{\mathsf {N}})\gamma _s(C_{\mathsf {N}}).
\end {gather*}
Taking $s=1=t$ in the above equalities we get
$C_{\mathsf {N}}^2=C_{\mathsf {N}}$.

To prove (3) note first that Proposition~\ref {prop:pparest},
  together with (1) and (2), show that $\alpha $, $\alpha '$ and
  $\theta $ are product partial actions. The domains of $\alpha $, $\alpha '$ and
  $\theta $ are respectively $D_t=A\beta _t(A)$, $D_t'=A'\beta _t'(A')$,
  and $\bar {D}_t=C_{\mathsf {N}}\gamma _t(C_{\mathsf {N}})
  =\begin {pmatrix}A\beta _t(A)&A\beta _t(A)X\\
  A'\beta _t'(A')Y&A'\beta _t'(A')\end {pmatrix}$,
 $\forall t\in G$. So it is clear that
  $\alpha =\theta |_A$ and $\alpha '=\theta |_{A'}$. Finally, by (1):
\[ Y_1D_tX_1
=YA\beta _t(A)AX
=A'\beta _t'(A')
=D'_t,\ \  \forall t\in G,\]
which ends the proof.
\end {proof}

\begin {thm}\label {thm:strongmoritaE}
Let $\mathsf {M}=(A, B, X, Y, \tau _A, \tau _B)$
be a strong-graded-equivalence between partially-strongly-$G$-graded
algebras $A$ and
$B$, and let $\gamma ^A$ and $\gamma ^B$ be the
canonical partial actions of $G$ on $\ps {A}$ and $\ps {B}$ respectively. Then
$\gamma ^A$ and $\gamma ^B$ are Morita equivalent partial actions.
\end {thm}
\begin {proof}
Let $C=\oplus _{t\in G}C_t$ be the graded Morita ring of
$\mathsf {M}$. By  Theorem~\ref {thm:moritaE}  we know that
$(\mathsf {M}\# G,\beta ^C)$ is a Morita equivalence between $\beta ^A$
and $\beta ^B$. The ideal $\ps {A}$ of $A\# G$ is idempotent and
$\ps {A}\beta ^A_t(\ps {A})=\beta ^A_t(\ps {A})\ps {A}$, $\forall t\in G$
(recall Proposition~\ref {prop:restsmash}). Thus all we need to do is
to show that $(Y\# G)\ps {A}(X\# G)=\ps {B}$, and then use
Proposition~\ref {prop:reseqppa}. By Proposition~\ref {prop:strong} and
Corollary~\ref {cor:sgt-sgeq}, we have for each $r,s\in G$:
\begin {gather*}
\sum _{u,v\in G}(Y_{r\m
  u}\munit _{r,u})(A_{u\m }A_v\munit _{u,v})(X_{v\m s}\munit _{v,s})
=\sum _{u,v\in G}Y_{r\m u}A_{u\m }A_vX_{v\m s}\munit _{r,s}\\
=\sum _{v\in G}Y_{r\m }A_vX_{v\m s}\munit _{r,s}
=Y_{r\m }X_s\munit _{r,s}=B_{r\m }B_s.
\end {gather*}
Then $(Y\# G)\ps {A} (X\# G)=\ps {B}$.

\begin {comment}
We compute
$C_{r\m }C_s$:
\begin {gather*}
C_{r\m }C_s
%=\begin {pmatrix} A_{r\m }  & X_{r\m } \\ Y_{r\m }& B_{r\m } \end {pmatrix}
%\begin {pmatrix} A_{s}  & X_{s} \\ Y_{s}& B_{s} \end {pmatrix}\\
=\begin {pmatrix} A_{r\m }A_s +\tau _A(X_{r\m }\otimes _B Y_s) &
  A_{r\m }X_s+X_{r\m }B_s \\ Y_{r\m }A_s+B_{r\m }Y_s&
  \tau _B(Y_{r\m }\otimes _AX_s)+B_{r\m }B_s \end {pmatrix}.
\end {gather*}
Then, using Proposition~\ref {prop:strong} to simplify the latter
matrix we get
\[C_{r\m }C_s=\begin {pmatrix} A_{r\m }A_s &
  A_{r\m }A_sX_1 \\ B_{r\m }B_sY_1&
  B_{r\m }B_s \end {pmatrix}
=\begin {pmatrix} A_{r\m }A_s &
  X_1B_{r\m }B_s \\ Y_1A_{r\m }A_s&
  B_{r\m }B_s \end {pmatrix}.\]
Thus
$\oplus _{r,s}C_{r\m }C_s\munit _{r,s}=\begin {pmatrix}\oplus _{r,s} A_{r\m }A_s\munit _{r,s} &
  \oplus _{r,s}A_{r\m }A_sX_1\munit _{r,s} \\ \oplus _{r,s}B_{r\m }B_sY_1\munit _{r,s}&
  \oplus _{r,s}B_{r\m }B_s\munit _{r,s} \end {pmatrix}
,$
that is \[\ps {C}=\begin {pmatrix} \ps {A} &
  \ps {A}X_1 \\ \ps {B}Y_1&
  \ps {B} \end {pmatrix}
=\begin {pmatrix} \ps {A} &
  X_1\ps {B} \\ Y_1\ps {A}&
  \ps {B} \end {pmatrix}
.\]
}
\fac {
To simplify our notation we set $I:=\ps {C}$ and $\beta :=\beta ^C$. Also,
we will write
\[\beta _t\begin {pmatrix} a_{\#}  &
  x_{\#} \\ y_{\#}&
  b_{\#} \end {pmatrix}
=\begin {pmatrix} \beta _t(a_{\#})  & \beta _t(x_{\#})
  \\ \beta _t(y_{\#})& \beta _t(b_{\#}) \end {pmatrix}, \textrm { for all }\begin {pmatrix} a_{\#}  &
  x_{\#} \\ y_{\#}&
  b_{\#} \end {pmatrix}\in C\# G.\]
Note that
$\beta _t(\ps {A})X_1=\beta _t(\ps {A}X_1)=\beta _t(X_1\ps {B})=X_1\beta _t(\ps {B})$,
and
$\beta _t(\ps {B})Y_1=\beta _t(\ps {B}Y_1)=\beta _t(Y_1\ps {A})=Y_1\beta _t(\ps {A})$.
We calculate $I_t=I\beta _t(I)$:
\begin {gather*}
I_t=\begin {pmatrix} \ps {A} &
  \ps {A}X_1 \\ \ps {B}Y_1&
  \ps {B} \end {pmatrix}\begin {pmatrix} \beta _t(\ps {A}) &
  \beta _t(\ps {A}X_1) \\ \beta _t(\ps {B}Y_1)&
  \beta _t(\ps {B}) \end {pmatrix}
=\begin {pmatrix} I_t^A &
  I_t^AX_1 \\ \ps {B}_tY_1&
  I_t^B \end {pmatrix}
\end {gather*}
It is easy to check that $(I_t^A,I_t^B,I_t^AX_1,I_t^BY_1)$ is a Morita
equivalence. Moreover $\gamma ^C=\{\gamma _t:I_{\m }\to I_t\}_{t\in G}$ is a partial
action of $G$ on $C\# G$ that leaves $A\#G$, $B\#G$, $X\#G$ and $Y\#G$
invariant, and that restricted to $A\#G$ and $B\#G$ agrees with
$\gamma ^A$ and $\gamma ^B$ respectively, so the latter partial actions
are Morita equivalent.
\end {comment}
\end {proof}
{
\begin {cor}\label {cor:strongmoritaE}
Let $A$ and $B$ be partially-strongly-$G$-graded algebras, and let
$\gamma ^A$ and $\gamma ^B$ be their corresponding canonical partial
actions. Then $A$ and $B$ are strongly-graded-equivalent if and only
if $\gamma ^A$ and $\gamma ^B$ are Morita equivalent partial actions.
\end {cor}
\begin {proof}
Just combine Propositions~\ref {prop:sge is eq},~\ref {prop:moritaeqs} with
Theorems~\ref {thm:sg},~\ref {thm:strongmoritaE}.
\end {proof}
}
\subsection {{Weak equivalence}}\label {weakEqParAc}
{Proposition~\ref {prop:moritaeqs} suggests the following notion:
\begin {defn}\label {defn:weakeq}
Let $\alpha $ and $\alpha '$ be product partial actions on the algebras
$A$ and $A'$ respectively. We say that $\alpha $ and $\alpha '$ are
weakly equivalent whenever $A\rtimes _\alpha G$ and
$A'\rtimes _{\alpha '} G$ are graded-equivalent.
\end {defn}

 Of course Morita equivalence implies weak
equivalence. From Proposition~\ref {prop:ge is eq} we immediately obtain:
\begin {prop}\label {prop:we is eq} Weak equivalence of product partial
  actions is an equivalence relation.
\end {prop}
\begin {prop}\label {prop:EqStronGR}
Let $B$ and $B'$ be strongly-graded algebras. Then $B$ and $B'$ are
strongly-graded-equivalent if and only if they are graded-equivalent.
\end {prop}
\begin {proof}
We only need to prove the converse. Since $B$ and $B'$ are strongly-graded algebras we have $\ps {B}=B\# G$ and $I^{B'}=B'\# G$. Hence, by
Theorem~\ref {thm:sg}, $B$ is strongly-graded-equivalent to $(B\#
G)\rtimes _{\beta ^B}G$ and $B'$ is strongly-graded-equivalent to $(B'\#
G)\rtimes _{\beta ^{B'}}G$.
Furthermore, thanks to  (3)  of Theorem~\ref {thm:moritaE}
$(B\#G)\rtimes _{\beta ^{B}}G$ and $(B'\# G)\rtimes _{\beta ^{B'}}G$ are
strongly-graded-equivalent. Then, in view of
Proposition~\ref {prop:sge is eq}, we conclude that $B$ and $B'$ are
strongly-graded-equivalent.
\end {proof}

\section {Globalization of product partial actions}\label {sec:glob}
%\vspace *{1cm}
\par As we show in the present section, much of the work done for
  regular partial actions in
\cite [Sections 4 and 5]{ades} can be extended to the case of product
partial actions.

{\begin {defn}\label {defn:meglob}
      Let $\alpha $ be a product partial action, and $\beta $ a global
      action. We say that $\beta $ is a Morita enveloping action of
      $\alpha $ if $\beta $ is a minimal globalization of a product
      partial action $\alpha '$ which is Morita equivalent with
      $\alpha .$
      \end {defn}
}

\begin {thm}\label {thm:moritaglob} Let $\alpha = \{
  {D}_{t\m }\stackrel {{\alpha }_t}{\to }  {D}_t\}_{t \in G}$
  be a product partial action of a group $G$ on an algebra $A .$
Then $\alpha $ {has a Morita enveloping action.}
More precisely, if $B=A\rtimes _\alpha G$, and
  $\alpha ':=\gamma ^{B}$ is
the canonical partial action of $G$ on $\ps {B}$, and
{$\beta ':=\beta ^{B}$ is the canonical action on $B\#G$},
then $\alpha $ and
$\alpha '$ are Morita equivalent product partial actions, {and
  $\beta '$ is a Morita enveloping action for $\alpha $.}
 \end {thm}
 \begin {proof}
 {Note that $\beta '$ is a minimal
globalization of $\alpha '=\beta '|_{\ps {B}}$ by Proposition~\ref {prop:restsmash}}. Now it is very easy to adapt the proof of
\cite [Theorem~4.1]{ades}, keeping in mind  Remark~\ref {def2MoritaEq}. Indeed, obviously, there is no need to prove
property (2) from
\cite {ades} for the ideals $D'_t= \ps {B} \beta _t^B (\ps {B}) .$ Nevertheless, equality  (12) from \cite {ades} should be used in the proof.
 However,   the latter  is an immediate consequence of equality (9)
 from \cite {ades}, which in our case is given by
 Proposition~\ref {prop:restsmash}. The rest of the proof goes without
 any change.
\end {proof}

In what follows, given a product partial action $\alpha $ on an
  algebra $A$, we
  will denote by $\beta ^{\alpha }$ the canonical action of $G$ on
  $(A\rtimes _\alpha G)\#G$, and by
  $\gamma ^\alpha =\beta ^{\alpha }|_{\ps {(A\rtimes _\alpha G)}}$ the
  canonical partial action of $G$ on $\ps {(A\rtimes _\alpha G)}$. Thus,
according with Theorem~\ref {thm:moritaglob},
$\alpha \stackrel {M}{\sim }\gamma ^\alpha $, and $\beta ^\alpha $ is a
Morita enveloping action for $\alpha $. We refer to $\beta ^\alpha $ as
the \textit {canonical Morita enveloping action of} $\alpha $. Note
that, in virtue of the comments preceding Remark~\ref {rem:eneveloping ac}
%Proposition~\ref {prop:pparest}
and
Proposition~\ref {prop:duality}, the correspondences that send
$\alpha $ to $\gamma ^\alpha $ and to $\beta ^\alpha $, as well as to the
corresponding skew group algebras, determine  functors.

\par As announced before, the converse of Proposition~\ref {prop:moritaeqs} holds:
\begin {thm}\label {thm:moritaskew} Let $\alpha = \{
  {D}_{t\m }\stackrel {{\alpha }_t}{\to }  {D}_t\}_{t \in G}$ and
  $\alpha ' = \{
  {D'}_{t\m }\stackrel {{\alpha '}_t}{\to }  {D'}_t\}_{t \in G}$ be
  product  partial actions of $G$ on  algebras $A$ and  ${A
  }',$ respectively. Then the skew group algebras $A \rtimes _{\alpha }
  G$ and ${A}' \rtimes _{{\alpha }'} G$ are strongly-graded-equivalent
  if and only if $\alpha $ and $\alpha '$ are Morita equivalent.
\end {thm}
\begin {proof}
  For the {\it `only if'} part, we have
  $\alpha \stackrel {M}{\sim }\gamma ^\alpha $ and
  $\alpha '\stackrel {M}{\sim }\gamma ^{\alpha '}$ by
  Theorem~\ref {thm:moritaglob} and
  $\gamma ^\alpha \stackrel {M}{\sim }\gamma ^{\alpha '}$ by
  Theorem~\ref {thm:strongmoritaE}. Then
  $\alpha \stackrel {M}{\sim }\alpha '$ by
  Proposition~\ref {prop:equivalencerel}.
  The {\it `if'} part  is
	Proposition~\ref {prop:moritaeqs}.
\end {proof}

\begin {cor}\label {cor:wthenM}
Two global actions on idempotent algebras are weakly equivalent if and only if they are Morita equivalent.
\end {cor}
\begin {proof} The claim follows at once from Proposition~\ref {prop:EqStronGR} and Theorem~\ref {thm:moritaskew} above.
\end {proof}

\par {The next two results will show that any Morita enveloping
  action of a product
  partial action $\alpha $ is Morita equivalent to the canonical
Morita enveloping action $\beta ^\alpha $ of $\alpha $.}

\begin {thm}\label {thm:moritaenveloping}
Let $\alpha = \{ {\alpha }_t : {D }_{t\m } \to {D }_t\}_{t \in G}$ and
$\alpha ' = \{ {\alpha '}_t : {D '}_{t\m } \to
{D }'_t\}_{t \in G}$ be Morita equivalent product  partial actions
of $G$ on algebras $A $ and  $A '$, respectively.
Then $\beta ^\alpha \stackrel {M}{\sim }\beta ^{\alpha '}$.
\end {thm}
\begin {proof}
By Proposition~\ref {prop:moritaeqs} the skew group
 algebras $A \rtimes _{\alpha } G$ and ${A}'
\rtimes _{{\alpha }'} G$ are strongly-graded-equivalent. Then our claim
follows from (3) of Theorem~\ref {thm:moritaE}.
\end {proof}
\begin {prop}\label {prop:envmorenv}
Let $\alpha = \{ {\alpha }_t : {D }_{t\m } \to {D }_t\}_{t \in G}$ be a product
partial action of $G$ on $A  .$
%that satisfies \eqref {D's}.
If $\beta :G\times B\to B$ is a minimal globalization of $\alpha $, {then
$\beta \stackrel {M}{\sim }\beta ^\alpha $.}
\end {prop}
\begin {proof}
The skew group algebras $A\rtimes _\alpha G$ and $B\rtimes _\beta G$ are
graded-equivalent by Theorem~\ref {thm:crossed product ppa and
  env}. Then {$\beta ^{\alpha }$ and $\beta ^{\beta }$} are Morita
equivalent by (3) of Theorem~\ref {thm:moritaE}. Since $B\rtimes _\beta
G$ is strongly-graded, and $I^{B\rtimes _\beta G}=(B\rtimes _\beta G)\#
G$, then $B\rtimes _\beta G$ is in fact strongly-graded-equivalent to
$\big ((B\rtimes _\beta G)\# G\big )\rtimes _{{\beta ^{\beta }}}G$ by
Corollary~\ref {cor:sg}. Thus $\beta $ and {$\beta ^{\beta }$} are Morita
equivalent by Theorem~\ref {thm:moritaskew}. Since Morita equivalence
of actions is an equivalence relation (recall
Proposition~\ref {prop:equivalencerel}), it follows that
{$\beta ^{\alpha }$} and $\beta $ are Morita equivalent.
\end {proof}
\begin {cor}\label {cor:umg}
    If $\beta $ and $\beta '$ are minimal globalizations of the Morita
    equivalent product partial actions $\alpha $ and $\alpha '$
    respectively, then $\beta $ and $\beta '$ are Morita equivalent
    actions.
\end {cor}
\begin {proof}
This follows at once from Theorem~\ref {thm:moritaenveloping} and
Proposition~\ref {prop:envmorenv}.
\end {proof}

%\sout {The proof of }\cite [Proposition 5.5]{ades} \sout {readily works for product partial actions and
%shows that $\tilde {\alpha }$ is indeed a Morita enveloping action for $\alpha .$
%With the above definition, }
Our previous results can be summarized as
follows:
\begin {thm}\label {thm:main2}
Let $\alpha = \{ {\alpha }_t : {D }_{t\m } \to {D }_t\}_{t \in G}$ be a product
partial action of $G$ on $A  .$
%that satisfies \eqref {D's}.
Then $\alpha $
has a Morita enveloping action, which is unique up to Morita
equivalence. Moreover, for every Morita enveloping action
$\beta :G\times B \to B $ of $\alpha $, the skew group algebras $A \rtimes _\alpha
G$ and $B \rtimes _\beta G$ are graded-equivalent.
\end {thm}
\begin {proof}
Theorem~\ref {thm:moritaglob} ensures the existence of a Morita
enveloping action for $\alpha $, and its uniqueness up to Morita
equivalence follows from Proposition~\ref {prop:equivalencerel} and
Corollary~\ref {cor:umg}.
%is given by Theorem~\ref {thm:moritaenveloping}.
The last
claim is a consequence of Theorem~\ref {thm:crossed product ppa and
  env}, Proposition~\ref {prop:moritaeqs} and Proposition~\ref {prop:ge
  is eq}.
%Proposition~\ref {prop:envmorenv}.
\end {proof}

{To conclude the section, we summarize several characterizations of
  Morita and weak equivalences of product partial actions.
\begin {prop}\label {prop:charmoreq}
Let $\alpha $ and $\alpha '$ be product partial actions of $G$ on the
algebras $A$ and $A'$ respectively. Then the following are equivalent:
\begin {enumerate}
 \item $\alpha $ and $\alpha '$ are Morita equivalent.
 \item $A\rtimes _\alpha G$ and $A'\rtimes _{\alpha '} G$ are strongly-graded-equivalent.
 \item $\gamma ^\alpha $ and $\gamma ^{\alpha '}$ are Morita equivalent.
 \item $\ps {A}\rtimes _{\gamma ^\alpha }G$ and
   $\ps {A'}\rtimes _{\gamma ^{\alpha '}}G$ are strongly-graded-equivalent.
\end {enumerate}
\end {prop}
\begin {proof}
Theorem~\ref {thm:moritaskew} implies that (1) and (2) are equivalent,
as well as (3) and (4). Finally, by Corollary~\ref {cor:strongmoritaE}, (2) and (3) are equivalent.
\end {proof}
\begin {prop}\label {prop:charweakeq}
Let $\alpha $ and $\alpha '$ be product partial actions of $G$ on the
algebras $A$ and $A'$ respectively, with corresponding Morita
enveloping actions $\beta $ and $\beta '$, acting on $B$ and
$B'$ respectively. Then the following are equivalent:
\begin {enumerate}
 \item $\alpha $ and $\alpha '$ are weakly equivalent.
 \item $\beta $ and $\beta '$ are weakly equivalent.
 \item $\gamma ^\alpha $ and $\gamma ^{\alpha '}$ are weakly equivalent.
 \item $\beta ^\alpha $ and $\beta ^{\alpha '}$ are weakly equivalent.
 \item $\alpha $ and $\beta '$ are weakly equivalent.
 \item $A\rtimes _\alpha G$ and $A'\rtimes _{\alpha '} G$ are graded-equivalent.
 \item $A\rtimes _\beta G$ and $A'\rtimes _{\beta '} G$ are graded-equivalent.
 \item $\ps {A}\rtimes _{\gamma ^\alpha }G$ and
   $\ps {A'}\rtimes _{\gamma ^{\alpha '}}G$ are graded-equivalent.
 \item $((A\rtimes _{\alpha }G)\#G)\rtimes _{\beta ^{\alpha }} G$ and
   $((A'\rtimes _{\alpha '}G)\#G)\rtimes _{\beta ^{\alpha '}} G$ are graded-equivalent.
 \item $A\rtimes _\alpha G$ and $B'\rtimes _{\beta '} G$ are graded-equivalent.
\end {enumerate}
Moreover, by Corollary~\ref {cor:wthenM}, (2) and (4) above can be
replaced by (2') and (4') below:
\begin {enumerate}
 \item [(2')] $\beta $ and $\beta '$ are Morita equivalent.
 \item [(4')] $\beta ^\alpha $ and $\beta ^{\alpha '}$ are Morita equivalent.
\end {enumerate}
\end {prop}
\begin {proof}
 Combining Theorems~\ref {thm:crossed product
  ppa and env}, \ref {thm:moritaglob}, \ref {thm:main2} and
Proposition~\ref {prop:we is eq},
and recalling that Morita equivalence implies weak equivalence, we
see that
$\alpha $, $\beta $, $\gamma ^\alpha $ and  $\beta ^\alpha $ are weak
equivalent, as well as
$\alpha '$, $\beta '$, $\gamma ^{\alpha '}$, $\beta ^{\alpha '}$. Hence the
first five assertions are equivalent to each other. Finally, according to the definition of weak equivalence, the last five
sentences are just rephrasings of the first five ones, all the ten
sentences are equivalent.
\end {proof}
}
%\vspace *{1cm}
%\textcolor {red}{\rule {\textwidth }{10pt}}

\section {Stabilization of graded algebras.}\label {sec:stab}

%We shall show in this section that stronlgy graded-equivalent partially-strongly-graded algebras with ortogonal local units   become graded isomorphic after stabilization.
Following  \cite {des}  we shall  say that a (non-necessarily graded)
algebra $A$ possesses {\it orthogonal local units} if there exists a set
of (non-necessarily central) pairwise orthogonal idempotents $E$
in $A$ such that

\begin {equation}\label {ortolocal1}
 A = \bigoplus _{e\in E}\; A e = \bigoplus _{e\in E} \; e A .
\end {equation}

\noindent Algebras  $A $ with (\ref {ortolocal1}) are also called
algebras with enough idempotents (see \cite {Fuller}). Note that algebras with  orthogonal local units generalize
algebras with a countable set of local units (see   \cite [p. 3300]{des}).
%It is also known by a result
%of P. N. \'{A}nh and L. M\'{a}rki \cite {AnhMarki} that every Morita
%equivalence class of rings with local units contains rings with orthogonal local units.

It is proved in
   \cite {des} that if $A$ and $B$ are Morita equivalent algebras with  orthogonal local units, then there is an
isomorphism of algebras
\begin {equation}\label {stableiso}
{\rm FMat}_{\mathcal X } (A) \cong {\rm FMat}_{\mathcal X} (B),
\end {equation} where ${\mathcal X}$ is an  appropriately chosen  infinite set.
Furthermore, it is shown in \cite {ades} that if  $A$ and $B$ are skew group algebras of Morita equi\-va\-lent
regular partial actions of a group $G$ on algebras with  orthogonal local units, then the isomorphism in
\eqref {stableiso} is graded.  In the latter case,  we now  know  by
Theorem~\ref {thm:moritaskew} that $A$ and $B$ are strongly-graded-equivalent.

Given a $G$-graded algebra $A$ and a set of indexes $\mathcal X,$ define a  $G$-grading on
${\rm FMat} _{\mathcal X} (A) $ by ta\-king the  $g$-homogeneous component of $ {\rm FMat} _{\mathcal X} (A) $ to be
${\rm FMat} _{\mathcal X}({A}_g),$ $g \in G.$  Note that this grading is different from that one considered in Theorem~\ref{thm:duality}.
We give the next:

\begin {thm}\label {thm:gradedstableiso}
Let
$$A = \bigoplus _{e\in E}\; A e = \bigoplus _{e\in E} \; e A ,$$ and
$$ B = \bigoplus _{f\in F}\; B f = \bigoplus _{f\in F} \; f B,$$  be partially-strongly-$G$-graded algebras with
  orthogonal local units.  Suppose that $A$ and $B$ are  strongly-graded-equivalent. Then  for any
infinite  set of indexes $\mathcal X,$ whose cardinality is bigger  than or
equal to those of $E$ and $F,$ there exists a graded isomorphism of algebras
\eqref {stableiso}.
% \begin {equation} \label {eq:gradedstableiso}
 %{\rm FMat} _{\mathcal X} (A) \cong {\rm FMat} _{\mathcal X}(B).
%\end {equation}
\end {thm}

\begin {proof}
Let $(A,B, {}_A X _B , {}_B Y_A,\tau _A,\tau _B)$ be a
  strong-graded-equivalence between $A$ and $B.$ In particular, $A$ and $B$ are Morita equivalent as non-necessarily graded algebras,
so that
\cite [Corollary 8.4]{des} implies the existence of an isomorphism of algebras
\eqref {stableiso}, and we need to check that it is a graded isomorphism. If
$A$ and $B$ are skew group algebras of  regular partial actions of  $G$ on algebras with orthogonal local units, this was
verified in \cite [Theorem 6.1]{ades}   by showing that the maps involved in  the construction of  \eqref {stableiso} are all
graded. It turns out that the arguments given in the proof of  \cite [Theorem 6.1]{ades} can be easily adapted to our
more general case. The adaptation has to be done at the starting point of the process and the subsequent steps follow
the same way.
 Since we have a strong-graded-equivalence between $A$ and $B, $ and these algebras are partially-strongly-graded, it
follows by Corollary~\ref {cor:sgt-sgeq} that  $(A_1,B_1,X_1,Y_1,\tau _A^1,\tau _B^1)$ is a Morita equivalence.
Obviously,  $E\subseteq A_1$ and $F \subseteq B_1.$ Consequently, for any $f\in F,$ using the trace map  $\tau _B^1,$
we may write
 $f =\sum _{i=1}^{n_f} y_i   x_i,$ where $x_i =
x^{(f)}_i \in X_1,$ and $ y_i = y^{(f)}_i  \in Y_1,$ $x_i = x_i f$, and $f y_i =y_i$ for all $i.$
Then it is readily seen that the map  $${\pi }_f:{A}^{n_f} \ni (r_1, r_2,  \ldots , r_{n_f}) \mapsto
\sum r_i x_i   \in X  f $$    is a graded epimorphism of left ${A}$-modules, as well as its splitting map
$${\rho }_f : X f \ni y f  \mapsto ( y f y_1 , y f y_2  , \ldots , y f y_{n_f}   ) \in
{A}^{n_f}.$$  Taking $K_f = {\rm Ker}\,  {\pi }_f$ and denoting by ${\mu
}_f$  the  embedding of $K_f$ in $  {A}^{n_f} ,$ we see that  all maps in the exact   sequences
\begin {equation*}\label {exact1}
0 \longrightarrow K_{f} \stackrel {{\mu }_f}\longrightarrow
{A}^{n_f} \stackrel {{\pi }_f}\longrightarrow X f  \longrightarrow
0,
\end {equation*}  and
\begin {equation*}\label {exact2}
0 \longleftarrow K_{f} \stackrel {{\tau }_f}\longleftarrow
{A}^{n_f} \stackrel {{\rho }_f}\longleftarrow X f   \longleftarrow 0,
\end {equation*}  preserve the $G$-gradings.
Similar maps are constructed replacing $A$ by $B$ and $X$ by $Y.$ The subsequent steps involve the use of the
functor $X_{B} \otimes _{-}$, application of the above maps to direct summands, rearrangements
of direct summands, each time resulting in graded maps, and leading  to a graded isomorphism
$A^{(\mathcal X)} \to X^{(\mathcal X)}$ of left $A$-modules, which is finitely determined in the sense of
\cite [Definition 7.3]{des}, and whose inverse is also finitely determined. Here
$A^{(\mathcal X)}$ (respectively,  $X^{(\mathcal X)}$) stands for the direct sum of copies of $A$ (respectively,  $X$),
indexed by the elements of $\mathcal X.$
An important point is to interpret the finitely determined isomorphism $A^{(\mathcal X)} \to X^{(\mathcal X)}$ as a
so-called   row and column summable ${\mathcal X} \times {\mathcal X}$-matrix   $[\psi ]$
over   $ {\rm RCFMat}_{(e,f) \in E \times F} (e X f ):$
$$[{\psi }] \in {\rm RCSumMat}_{X}({\rm RCFMat}_{(e,f) \in E \times F} (e X f )).$$ Then $[{\psi }]$ is used to
define  the  maps $$\Psi : {\rm
FMat}_{\mathcal X}( A ) \to {\rm FMat}_{\mathcal X}( X ), \; \; \; \mbox {and} \; \;
\; {\Psi }' : {\rm FMat}_{\mathcal X} ( B ) \to {\rm FMat}_{\mathcal X} ( X ), $$
and, with the help of some  preliminary results from  \cite {des} the isomorphism
\eqref {stableiso} is obtained as the composition $({\Psi }') \m \circ {\Psi }.$ The fact that
$A^{(\mathcal X)} \to X^{(\mathcal X)}$ is graded implies that each entry of  $[{\psi }]$ belongs to the
$1$-homogeneous component
$e X_1 f$ of $e X f.$ The latter yields that
\eqref {stableiso} is graded. The details can be seen in the proofs of \cite [Theorem 8.2, Corollary 8.4]{des} and
in the comments to them given to justify \cite [Theorem 6.1]{ades}. \end {proof}

\section*{Acknowledgments}
 The authors thank the  referee for useful comments.
F.A. was partially supported by Fapesp grant 2017/23242-0.
 M.D. was partially supported by Fapesp grants  2015/09162-9, 2020/16594-0 and by CNPq grant 312683/2021-9.
        R.E. was partially supported by Fapesp grant 2017/26645-9 and by CNPq grant 303173/2020-3.


\begin {thebibliography}{10}

 \bibitem {AlbuNast} T.\ Albu, C.  N\u {a}st\u {a}sescu,
  \textit {Infinite group-graded rings, rings of endomorphisms, and
    localization,} J. Pure Appl. Algebra \textbf {59}  (1989), no. 2,
  125--150.

 \bibitem {aa} B. Abadie and F. Abadie, \textit {Ideals in
     cross sectional C*-algebras of Fell bundles}, Rocky Mountain
   J. Math. 47 (2017), no. 2, 351--381.

 \bibitem {abenv} F. Abadie, \textit {Enveloping actions and Takai
     duality for partial actions}, J. Funct. Anal. \textbf {197}
   (2003), 14--67.

 \bibitem {abf} F. Abadie, A.~Buss, D.~Ferraro, \textit {Morita
     enveloping Fell bundles}, Bull. Braz. Math. Soc.,
     New Series \textbf {50} (1), (2019) 3--35.

 \bibitem {af} F. Abadie, D.~Ferraro, \textit {Equivalence of Fell
     bundles over groups},  J. Operator Theory \textbf {81} (2), (2019)
   273--319.

 \bibitem {ades} F.~Abadie, M.~Dokuchaev, R.~Exel, and J.~J.~Sim\'on,
   \textit {Morita equivalence of partial group actions and
     globalization},    Trans. Amer. Math. Soc. \textbf {368} (2016),
   4957--4992.

\bibitem {Abrams} G.\ D.\  Abrams, \textit {Morita equivalence for rings with local units}, Commun. Algebra
\textbf {11} (1983), 801--837.

\bibitem {AnhMarki} P.\ N.\ \'{A}nh, L.\ M\'{a}rki, \textit {Morita equivalence
for rings without identity,}  Tsukuba J. Math. {\bf 11}
(1987), 1--16.



\bibitem{Bavula2009} V.\ V.\ Bavula, \textit{The Jacobian algebras,} J. Pure Appl. Algebra {\bf 213} (2009), 664--685.

\bibitem{Bavula2010} V.\ V.\ Bavula, \textit{The algebra of one-sided inverses of a polynomial algebra,} J. Pure Appl. Algebra {\bf 214}, (2010),  1874--1897.

\bibitem{Bavula2011} V.\ V.\ Bavula, \textit{The algebra of integro-differential operators on a polynomial algebra,}
J. Lond. Math. Soc., {\bf 83}, (2011), 517--543.



 \bibitem {beattiedual} M. Beattie, \textit {Duality Theorems for Rings
     with Actions or Coactions}, J.~Algebra \textbf {115} (1988),
     303--312.

\bibitem {Boisen} P. Boisen, \textit {Graded Morita Theory}, J.~Algebra
     \textbf {164} (1994), 1--25.

\bibitem {BussMeyerZhu} A. Buss, R. Meyer, and C. Zhu, \textit {A higher category approach to twisted
actions on C*-algebras},
  Proc. Edinb. Math. Soc. (2) \textbf{56} (2013), no. 2, 387--426.

\bibitem {cm} M. Cohen, S. Montgomery, \textit {Group-graded rings,
     smash products, and group actions},
   Trans. Amer. Math. Soc. \textbf {282} (1984), 237--258.

\bibitem {ADelRio91} A. del R\'{\i }o,
\textit {Graded rings and equivalences of categories},
Commun. Algebra \textbf {19}   (1991), no. 3, 997--1012.

\bibitem {ADelRio92} A. del  R\'{\i }o,
\textit {Categorical methods in graded ring theory},
Publ. Mat., Barc. \textbf {36},  (1992), no. 2A, 489--531.

 \bibitem {dokex} M. Dokuchaev, R. Exel, \textit {Associativity of
     crossed products by partial actions, enveloping actions and
     partial representations},  Trans. Amer. Math. Soc. \textbf {357}
     (2005), no. 5, 1931--1952.

\bibitem {DEP} M.\ Dokuchaev, R.\ Exel, P.\ Piccione,
                  \textit {Partial representations and partial group algebras},
                   J. Algebra, {text\bf 226} (2000), no. 1, 505--532.

 \bibitem {des} M. Dokuchaev, R. Exel, J.~J.~Sim\'on, \textit {Crossed
     products by twisted partial actions and graded algebras},  J. Algebra, {\textbf 320}   (2008),  no. 8,  3278--3310.

 \bibitem {DdRS} M. Dokuchaev, A. del R{\'\i }o, J. J. Sim\'{o}n, \textit {Globalizations of
partial actions on non unital rings},
 Proc. Am. Math. Soc., {\textbf 135} (2007),  no. 2,  343--352.

\bibitem {extwist} R. Exel, {\it Twisted partial actions: a
    classification of regular C*-algebraic bundles}, Proc. London
  Math. Soc. (3) \textbf {74} (1997), no. 2, 417--443.

 \bibitem {pdsfba} R.~Exel, \textit {Partial dynamical systems, Fell bundles and applications},
   Mathematical Surveys and Monographs, 224. American Mathematical Society, Providence, RI, 2017. vi+321 pp.

\bibitem {fd} J. M. Fell, R. S. Doran, \textit {``Representations of
    *-algebras, locally compact groups, and Banach *-algebraic
    bundles''}, Pure and Applied Mathematics vol. {\bf 125} and \textbf {126}, Academic Press, 1988.

\bibitem {Fuller} K.\ R.\ Fuller, \textit {On rings whose left modules are
direct sums of finitely generated modules},  Proc.\ Amer.\
Math.\ Soc, {\textbf 54} (1976), 39--44.

 \bibitem {garsim}J. L. Garc\'\i a, J. J. Sim\'on, \textit {Morita equivalence
    for idempotent rings}, J. Pure Appl. Algebra \textbf {76} (1991),
  no. 1, 39--56.

\bibitem {GordonGreen} R.\ Gordon,  E.\ Green,  \textit {Graded Artin algebras}, J. Algebra  \textbf {76}  (1982), 111--137.

\bibitem {Haefner94} J. Haefner,
\textit {Graded Morita theory for infinite groups},
J. Algebra, \textbf { 169}, (1994),   no. 2, 552--586.

\bibitem {Haefner95} J. Haefner, \textit {Graded equivalence theory with applications},
J. Algebra \textbf {172},  (1995), no. 2,   385--424.

\bibitem {MeniniNast} C.\ Menini and C.\ N\u {a}st\u {a}sescu, \textit {When is R-gr equivalent to the category of
modules?}, J. Pure  Appl. Algebra \textbf {51} (1988), 277--291.

\bibitem {nastoyst} C. N\u {a}st\u {a}sescu, F. Van Oystaeyen,
 ``Methods of graded rings.''
Lecture Notes in Mathematics 1836. Berlin: Springer,  (2004).

 \bibitem {nop} P.~Nystedt, J.~\"{O}inert and H.~Pinedo, \textit {Epsilon strongly graded rings, separability and semisimplicity},
 J. Algebra  \textbf {514}  (2018), 1--24.
% arXiv:1606.07592v1.

 \bibitem {quinn} D. Quinn, \textit {Group-graded rings and duality},
   Trans. Amer. Math. Soc. \textbf {292} (1985), 155--167.

 \bibitem {sehnem} C.~F.~Sehnem, \textit {Uma classifica\c c\~ao de fibrados de Fell est\'aveis},
   Master's thesis, Universidade Federal de Santa Catarina, 2014.


 \end {thebibliography}
\end {document}